\documentclass[reqno]{amsart}
\usepackage{amsfonts}
\usepackage{amsmath}
\usepackage{amsthm}
\usepackage{amssymb}
\usepackage{latexsym}

\newif\iffinal
\finaltrue
\newtheorem*{cor}{Corollary}
\newtheorem*{lem}{Lemma}
\newtheorem*{prop}{Proposition}

\theoremstyle{definition}
\newtheorem*{defn}{Definition}
\theoremstyle{definition}
\newtheorem{thm}{Theorem}
\newtheorem*{thm*}{Theorem}

\newtheorem*{rem}{Remark}

\newtheorem*{clm}{Claim}

\newenvironment{pf}{\proof}{\endproof}
\newenvironment{pf*}{\par\normalfont\topsep6pt plus 6pt\relax\trivlist%
\item[\hskip\labelsep\itshape\proofname.]\ignorespaces}{\endtrivlist}

\newcounter{cnt}
\newenvironment{enumerit}{\begin{list}{{\hfill\rm(\roman{cnt})\hfill}}{%
\settowidth{\labelwidth}{{\rm(iv)}}\leftmargin=\labelwidth%
\advance\leftmargin by 
\labelsep\rightmargin=0pt\usecounter{cnt}}}{\end{list}}

\theoremstyle{remark}

\numberwithin{equation}{section} 
\begin{document}

\newcommand{\thmref}[1]{Theorem~\ref{#1}}
\newcommand{\secref}[1]{Section~\ref{#1}}
\newcommand{\lemref}[1]{Lemma~\ref{#1}}
\newcommand{\propref}[1]{Proposition~\ref{#1}}
\newcommand{\corref}[1]{Corollary~\ref{#1}}
\newcommand{\remref}[1]{Remark~\ref{#1}}
\newcommand{\defref}[1]{Definition~\ref{#1}}
\newcommand{\er}[1]{(\ref{#1})}
\newcommand{\id}{\operatorname{id}}
\newcommand{\tensor}{\otimes}
\newcommand{\dist}{\operatorname{dist}}
\newcommand{\nc}{\newcommand}
\newcommand{\End}{\operatorname{End}}
\newcommand{\rnc}{\renewcommand}
\newcommand{\Maj}{\operatorname{Maj}}
\newcommand{\qbinom}[2]{\genfrac[]{0pt}0{#1}{#2}}
\nc{\cal}{\mathcal} \nc{\goth}{\mathfrak} \rnc{\bold}{\mathbf}
\renewcommand{\frak}{\mathfrak}
\renewcommand{\Bbb}{\mathbb}
\nc\bpi{{\mbox{\boldmath $\pi$}}}
\nc\bgamma{{\boldsymbol\gamma}}
\nc\wt{\operatorname{wt}}
\newcommand{\lie}[1]{\mathfrak{#1}}
\makeatletter
\def\section{\def\@secnumfont{\mdseries}\@startsection{section}{1}%
  \z@{.7\linespacing\@plus\linespacing}{.5\linespacing}%
  {\normalfont\scshape\centering}}
\def\subsection{\def\@secnumfont{\bfseries}\@startsection{subsection}{2}%
  {\parindent}{.5\linespacing\@plus.7\linespacing}{-.5em}%
  {\normalfont\bfseries}}
\makeatother
\def\subl#1{\subsection{\iffalse#1\fi}\label{#1}\def\subslbl{#1}}
\def\lbl#1{\label{\subslbl.#1}}
\def\loceqref#1{\eqref{\subslbl.#1}}
\nc{\wh}[1]{\widehat{#1}}
\nc{\wtil}[1]{\widetilde{#1}}
\nc{\bsm}[1]{\boldsymbol{#1}}

\nc{\Cal}{\cal} \nc{\Xp}[1]{X^+(#1)} \nc{\Xm}[1]{X^-(#1)}
\nc{\on}{\operatorname} \nc{\ch}{\mbox{ch}} \nc{\Z}{{\bold Z}}
\nc{\J}{{\cal J}} \nc{\C}{{\bold C}} \nc{\Q}{{\bold Q}}
\renewcommand{\P}{{\cal P}}
\nc{\N}{{\Bbb N}} \nc\boa{\bold a} \nc\bob{\bold b} \nc\boc{\bold
c} \nc\bod{\bold d} \nc\boe{\bold e} \nc\bof{\bold f}
\nc\bog{\bold g} \nc\boh{\bold h} \nc\boi{\bold i} \nc\boj{\bold
j} \nc\bok{\bold k} \nc\bol{\bold l} \nc\bom{\bold m}
\nc\bon{\bold n} \nc\boo{\bold o} \nc\bop{\bold p} \nc\boq{\bold
q} \nc\bor{\bold r} \nc\bos{\bold s} \nc\bou{\bold u}
\nc\bov{\bold v} \nc\bow{\bold w} \nc\boz{\bold z}

\nc\ba{\bold A} \nc\bb{\bold B} \nc\bc{\bold C} \nc\bd{\bold D}
\nc\be{\bold E} \nc\bg{\bold G} \nc\bh{\bold H} \nc\bi{\bold I}
\nc\bj{\bold J} \nc\bk{\bold K} \nc\bl{\bold L} \nc\bm{\bold M}
\nc\bn{\bold N} \nc\bo{\bold O} \nc\bp{\bold P} \nc\bq{\bold Q}
\nc\br{\bold R} \nc\bs{\bold S} \nc\bt{\bold T} \nc\bu{\bold U}
\nc\bv{\bold V} \nc\bw{\bold W} \nc\bz{\bold Z} \nc\bx{\bold X}

\nc\Ca{\cal A} \nc\Cb{\cal B} \nc\Cc{\cal C} \nc\Cd{\cal D}
\nc\Ce{\cal E} \nc\Cg{\cal G} \nc\Ch{\cal H} \nc\Ci{\cal I}
\nc\Cj{\cal J} \nc\Ck{\cal K} \nc\Cl{\cal L} \nc\Cm{\cal M}
\nc\Cn{\cal N} \nc\Co{\cal O} \nc\Cp{\cal P} \nc\Cq{\cal Q}
\nc\Cr{\cal R} \nc\Cs{\cal S} \nc\Ct{\cal T} \nc\Cu{\cal U}
\nc\Cv{\cal V} \nc\Cw{\cal W} \nc\Cz{\cal Z} \nc\Cx{\cal X}
\nc{\hyper}[2]{\genfrac{}{}{0pt}{}{#1}{#2}}
\def\<#1,#2>{#1(#2)}
\def\hatg{\wh{\lie g}}
\def\hath{\wh{\lie h}}
\def\hatbp{\wh{\mathbb P}}
\def\mbp{\mathbb P}
\newenvironment{enumspecm}[1]{\list{($#1\thecnt$)}{%
\settowidth{\labelwidth}{($#1{5}$)}%
\usecounter{cnt}\leftmargin=\labelwidth\advance\leftmargin by\labelsep%
\rightmargin=0pt\itemindent=0pt}}%
{\endlist}

\title[Path model for quantum loop modules]{Path model for quantum loop 
modules of fundamental type}
\iffinal
\author{Jacob Greenstein}
\thanks{The first author was supported by a 
Marie Curie Individual Fellowship of the European Community 
under contract no.~HPMF-CT-2001-01132}
\address{Institut de Math\'ematiques de Jussieu, Universit\'e
Pierre et Marie Curie, 175 rue du Chevaleret, Plateau 7D, F-75013
Paris, France}
\email{greenste@math.jussieu.fr}
\author{Polyxeni Lamprou}
\email{lamprou@edcsm.jussieu.fr}
\subjclass[2000]{Primary 17B67}
\fi
\date{\today}
\maketitle

\section{Introduction}

\subl{I10}
In this paper we construct a combinatorial realisation of
a certain class of simple integrable 
modules with finite dimensional weight spaces over a quantised affine algebra.

The best-known examples of such modules are the highest weight simple 
integrable modules~$V(\lambda)$. These modules are, essentially,
combinatorial objects for the following reasons. First of all, they 
can be defined for an arbitrary quantised Kac-Moody algebra. 
Next, the formal character of~$V(\lambda)$ is given by a universal formula
known as Kac-Weyl character formula (cf.~\cite[Chapter~10]{Kac}) and
determines~$V(\lambda)$ uniquely up to an isomorphism. 
Furthermore, $V(\lambda)$
is a quantum deformation (cf.~\cite{L1}) 
of a module over the corresponding Kac-Moody algebra which is also simple 
and has the same formal character. Finally, after~\cite{Ka1,L}, $V(\lambda)$
admits a crystal basis and a global basis. 

The properties of a crystal basis, formulated in an abstract way,
lead to the notion of a crystal as a set equipped with root operators
$\tilde e_\alpha$, $\tilde f_\alpha$ for each simple root~$\alpha$ of
the corresponding Kac-Moody algebra and some other operations which will
be discussed later. In particular, one associates with~$V(\lambda)$
a crystal~$B(\lambda)$ which encodes the major properties
of the module. For example, one can define, in a natural way, a tensor
product of crystals whose properties reflect these of the tensor product
of modules for the~$V(\lambda)$. Namely, a decomposition
of the tensor product of crystals~$B(\lambda)$ and $B(\mu)$ yields a
decomposition of~$V(\lambda)\tensor V(\mu)$.

\subl{I20}
The crystals~$B(\lambda)$ are known to admit numerous combinatorial
realisations. One of the most important, due to its simplicity and 
universality, is the path model of Littelmann~(cf.~\cite{Li94,Li95}).
In the framework of that model, $B(\lambda)$ is represented as a subset
of the set~$\mathbb P$ 
of piece-wise continuous linear paths in a rational vector subspace
a Cartan subalgebra of the Kac-Moody algebra connecting the origin
with an integral weight. Then the tensor product of crystals corresponds
to the concatenation of paths. The Isomorphism Theorem of Littelmann
(cf.~\cite{Li95}) stipulates that any subcrystal 
of~$\mathbb P$, which is generated over the associative monoid~$\Cm$ of 
root operators by a path which connects the origin with~$\lambda$ and
lies entirely in the dominant chamber, provides a realisation of~$B(\lambda)$.
Moreover, any two such realisations for~$\lambda$ fixed are isomorphic as 
crystals. In particular, they are isomorphic to the subcrystal of~$\mathbb P$
generated over~$\Cm$ by the linear path connecting the origin with~$\lambda$.

\subl{I30}
The case of affine Lie algebras is somewhat special since they admit,
besides the Kac-Moody presentation, an explicit realisation in
terms of loop algebras. 

Let~$\lie g$ be a finite dimensional simple Lie algebra of rank~$\ell$
over~$\bc$ and denote by~$\hatg $ the corresponding untwisted affine
algebra (cf.~\ref{P20}). 
Two quantum versions of~$\hatg $ are generally considered. They
will be denoted by~$\bu_q$ and~$\wh\bu_q$ respectively and differ
by a choice of the torus (cf.~\ref{P30}). The algebra $\bu_q$ can 
also be viewed as a subquotient of~$\wh\bu_q$. 

The algebra~$\bu_q$ admits finite dimensional integrable representations
which have been and still are being 
studied extensively (cf., to name but a few,
\cite{AK,CPqa,CP,Dr,EM,FM,FR,Ka2,N1,N2,VV}). These modules are parametrised
by $\ell$-tuples of polynomials over~$\bc(q)$ 
in one variable with constant term~$1$, known as Drinfel'd polynomials, and
are very different, in many respects, from highest weight integrable
modules. They are not, in general, determined by their formal character
(however, they are determined by their $q$-characters introduced 
in~\cite{FR}). They do not always admit classical limits and these limits,
when exist, are not necessarily simple modules over the corresponding 
affine Lie algebra and in fact may have a rather complex structure. 
Finally, it seems that existence of a crystal basis 
is an exception rather than a rule for this class of modules. The general
reason for these discrepancies is that the construction of finite 
dimensional~$\bu_q$~modules arises from the loop-like 
(Drinfel'd) presentation of~$\bu_q$ (cf.~\cite{Be,Dr,J}) peculiar 
to the Kac-Moody algebras of affine type.

\subl{I40}
Simple (infinite dimensional) integrable modules with finite dimensional
weight spaces were classified in~\cite{Cint,CPnew} for affine
Lie algebras and in~\cite{CG} for quantised affine algebras. 
Namely, such a module
is either a highest weight module~$V(\lambda)$ (or its graded dual) or
a loop module. The modules of the latter class are constructed,
in the quantum case, as simple submodules of the loop spaces of finite 
dimensional simple modules over~$\bu_q$. Namely, 
let~$\bpi=(\pi_1,\dots,\pi_\ell)$, $\pi_i\in\bc(q)[u]$
be an $\ell$-tuple of polynomials with constant term~$1$ and 
let~$V(\bpi)$ be the corresponding finite dimensional simple~$\bu_q$-module. 
Let~$m$ be the maximal positive integer 
such that all the~$\pi_i$, $i=1,\dots,\ell$ lie 
in~$\bc(q)[u^m]$. Then one can show (cf.~\cite{CG}) that the cyclic
group~$\bz/m\bz$ acts on the loop space 
$\wh V(\bpi):=V(\bpi)\tensor_{\bc(q)}\bc(q)[t,t^{-1}]$ and its action
commutes with that of~$\wh\bu_q$. In particular, simple submodules~$\wh
V(\bpi)^{(k)}$, $k=0,\dots,m-1$
correspond to distinct irreducible characters of the abelian group~$\bz/m\bz$.
We say that~$\wh V(\bpi)$ is of {\em fundamental type} 
if~$\pi_j(u)=\delta_{i,j}
(1-u^m)$ for some~$m>0$ and for some~$i\in\{1,\dots,\ell\}$ fixed. 
Henceforth we denote such an $\ell$-tuple of polynomials 
by~$\bsm\varpi_{i;m}$.

It turns out that simple submodules of~$\wh V(\bsm\varpi_{i;m})$ 
are determined by their formal characters up to a twist by an automorphism 
of~$\wh\bu_q$. In the present paper we show that these modules admit
a certain analogue of a crystal basis and construct a realisation in
the framework of Littelmann's path model
of the crystal associated to that basis in a natural way. 
The first example of~$\lie g$ of type~$A_\ell$, $m$ arbitrary and~$i=1$,
in which case the module~$V(\bsm\varpi_{i;1})$ is isomorphic to the
quantum analogue of the natural ($\ell+1$-dimensional) representation
of~$\lie g$ as a module over the quantised enveloping algebra~$U_q(\lie g)$
corresponding to~$\lie g$, was
considered in~\cite{G}. 
The case~$m=1$ was later 
treated, independently, by S.~Naito and D.~Sagaki (cf.~\cite{NS}) for~$\lie g$
of all types and for all~$i=1,\dots,\ell$. Here we consider all modules
of fundamental type for~$\lie g$ of all types which we believe 
to be the widest class of integrable $\wh\bu_q$-modules of level zero with
finite dimensional weight spaces which admit a combinatorial
realisation inside the path crystal of Littelmann. Our analysis is based
on the approach of~\cite{G} and on the results of~\cite{Ka2} and~\cite{NS}.

\subl{I50}
Let us briefly describe the principal results of this paper.
It was shown in~\cite{Ka2} that~$V(\bsm\varpi_{i;1})$ always
admits a crystal basis~$B(\bsm\varpi_{i;1})$ 
whose $m$th tensor power, for any~$m>0$
is indecomposable as a crystal. In order to treat the 
modules~$V(\bsm\varpi_{i;m})$
for an arbitrary~$m$ one has to 
introduce the notion of a $z$-crystal basis (cf.~\cite{CG} and~\defref{ZCB10}).
Roughly speaking, whilst crystal bases are preserved as sets by
the root operators of Kashiwara, 
$z$-crystal bases are preserved by these
operators up to a multiplication by a power of a complex number~$z$ 
only. Our first result is the following
\begin{thm}\label{thm1}
The simple module~$\wh V(\bsm\varpi_{i;m})^{(k)}$, 
$k=0,\dots,m-1$, $m>0$, admits
a $z$-crystal basis~$\wh B(\bsm\varpi_{i;m})^{(k)}$, 
where~$z$ is an $m$th primitive root of unity.
\end{thm}
From the combinatorial point of view multiplication of elements
of a basis by roots of unity is not important and one can get rid of it
associating a crystal to a $z$-crystal basis (cf.~\ref{ZCB12}).
It turns out that the crystal associated with~$\wh B(\bsm\varpi_{i;m})^{(k)}$ 
is indecomposable and
these are all indecomposable subcrystals of the affinisation (cf.~\ref{P50})
of the finite crystal~$B(\bsm\varpi_{i;1})^{\tensor m}$. That illustrates
once again how different loop modules are from highest weight modules.
Indeed, the affinisation of~$B(\bsm\varpi_{i;1})^{\tensor m}$ is also 
isomorphic
to the crystal basis of the simple $\wh\bu_q$-module 
$\wh V(\bpi)$ 
where~$\bpi=(\pi_1,\dots,\pi_\ell)$ with~$\pi_j(u)=\delta_{i,j}(1-u)^m$. 
Thus, the crystal basis of that simple module is a disjoint
union of indecomposable crystals.

Let~$\varpi_i$, $i=1,\dots,\ell$ be the fundamental
weights of~$\lie g$ extended by zero to weights of~$\hatg$ and
let~$\delta$ be the generator of imaginary roots of~$\hatg$ (cf.~\ref{P20}).
The main result of this paper is the following
\begin{thm}\label{thm2}
The associated crystal of~$\wh B(\bsm\varpi_{i;m})^{(k)}$ is isomorphic to
the subcrystal~$B(m\varpi_i+k\delta)$ of the Littelmann path crystal 
generated by the linear path connecting the origin with~$m\varpi_i+k\delta$.
\end{thm}

\subsection*{Acknowledgements}
We are greatly indebted to A.~Joseph who taught us all we know
about crystals. We are grateful to V.~Toledano-Laredo, and the
first author thanks B.~Leclerc, P.~Littelmann
and M.~Varagnolo, for numerous interesting discussions. 

\section{Preliminaries and notations}

\subl{P10}
Let~$\bc(q)$ be the field of rational functions in~$q$ with complex
coefficients, that is, the fraction field of~$\bc[q]$. Let~$\Ca\subset
\bc(q)$ be the ring~$\bc[q]$ localized at~$q=0$, which identifies
with the subring of rational functions in~$q$ regular
at~$q=0$.
Given~$m\ge n\ge0$, define
$$
[m]_q:=\frac{q^m-q^{-m}}{q^{}-q^{-1}},\quad
[m]_q!=[1]_q\cdots [m]_q,\quad
\qbinom{m}{n}_q:=\frac{[m]_q!}{[n]_q! [m-n]_q!}.
$$
All the above are Laurent polynomials in~$q$ over~$\bz$.

\subl{P20}
Set~$I=\{1,\dots,\ell\}$ and
let~$A=(a_{ij})_{i,j\in I}$ be the Cartan matrix of a finite dimensional
simple Lie algebra~$\lie g$ over~$\bc$ of rank~$\ell$. Fix a Cartan
subalgebra~$\lie h\subset\lie g$ and let~$\{\alpha_i\}_{i\in I}$
(respectively, $\{\alpha_i^\vee\}_{i\in I}$) be a basis of~$\lie h^*$
(respectively, of~$\lie h$) such that~$\<\alpha_i^\vee,\alpha_j>=
a_{ij}$. Define the fundamental weights~$\varpi_i\in\lie h^*$, 
$i\in I$ of~$\lie g$ by~$\<\alpha_i^\vee,\varpi_j>=\delta_{i,j}$,
where~$\delta_{i,j}$ is
the Kronecker's symbol, and 
let~$P_0$ be the free abelian group generated by the~$\varpi_i$, $i\in I$.
Let~$\theta=\sum_{i\in I} a_i\alpha_i$ be the highest root of~$\lie g$ with respect to~$\lie h$ and denote by~$\theta^\vee=\sum_{i\in I} a_i^\vee 
\alpha_i^\vee$ the corresponding co-root.

Set~$\wh I=I\cup\{0\}$ and let~$\wh A=(a_{ij})_{i,j\in\wh I}$ be
the generalised Cartan matrix of the untwisted affine Lie algebra~$\hatg $
associated with~$\lie g$. As a vector space,
$$
\hatg =\lie g\tensor_\bc\bc[t,t^{-1}]\oplus\bc c\oplus\bc\partial,
$$
where~$c$ is the canonical central element and~$\operatorname{ad}\partial=
t\frac{d}{dt}$. Then~$\hath=\lie h\oplus\bc c\oplus\bc\partial$
is a Cartan subalgebra of~$\hatg$. Set~$\alpha_0^\vee:=c-\theta^\vee$.
Define~$\delta\in\hath^*$ by~$\<\partial,\delta>=1$, $\delta|_{\lie h\oplus
\bc c}=0$ and set~$\alpha_0=\delta-\theta$. Then~$\{
\alpha_i\}_{i\in\wh I}$ (respectively, 
$\{\alpha_i^\vee\}_{i\in \wh I}$) is a set of simple roots
of~$\hatg$ and~$\<\alpha_i^\vee,\alpha_j>=a_{ij}$, $i,j\in\wh I$. 
Notice that~$\<\alpha_i^\vee,\delta>=0=\<c,\alpha_i>$
for all~$i\in\wh I$.

Define the fundamental weights~$\Lambda_i\in\hath^*$, $i\in\wh I$ 
of~$\hatg$ by conditions~$\<\alpha_i^\vee,\Lambda_j>=\delta_{i,j}$,
$\<\partial,\Lambda_i>=\delta_{i,0}$. Let~$P$ be the free
abelian group generated by the~$\Lambda_i$, $i\in\wh I$ and 
set~$\wh P:=P\oplus\bz\delta$. Extend the map
$\varpi_i\mapsto \Lambda_i-a_i^\vee \Lambda_0$ to an embedding 
of~$P_0$ into~$P$ and identify~$P_0$ with its image
inside~$P$ which in turn coincides with the set~$\{\lambda\in P\,:\,
\<c,\lambda>=0\}$. Let~$\xi:\wh P\to\wh P/\bz\delta$ be the canonical
projection. Notice that~$P$ identifies with~$\wh P/\bz\delta$ and
that~$\xi(\alpha_0)=-\theta$.

For all~$i\in\wh I$ define an elementary reflection~$s_i\in
\operatorname{Aut}\hath^*$
by~$s_i \lambda=\lambda-\<\alpha_i^\vee,\lambda>\alpha_i$ 
for all~$\lambda\in\hath^*$. The Weyl group $\wh W$ of~$\hatg$ (respectively,
the Weyl group~$W$ of $\lie g$) identifies
with the group generated by the~$s_i\,:\,i\in\wh I$ (respectively,
$i\in I$). The set of roots 
of~$\hatg$ is a disjoint union of the set of real roots~$\cup_{i\in
\wh I}\wh W\alpha_i$
and imaginary roots~$\bz\delta\setminus\{0\}$. If~$\beta$ is a real root,
denote the corresponding co-root by~$\beta^\vee$ and set~$s_\beta\lambda=
\lambda-\<\beta^\vee,\lambda>\beta$, $\lambda\in\hath^*$. Observe
that~$s_0=s_\theta$ as an automorphism of~$P$ and so~$\wh W$ identifies
with~$W$ when we consider the action of the former group on~$P$.

\subl{P30}
Let~$d_i$, $i\in\wh I$ be positive relatively prime integers such that
the matrix~$(d_i a_{ij})_{i,j\in\wh I}$ is symmetric and let~$q_i=
q^{d_i}$. Henceforth, for any symbol~$X_i$, $i\in\wh I$, set~$X_i^{(k)}:=
X_i^k/[k]_{q_i}!$.

The quantised affine algebra~$\wh\bu_q:=\bu_q(\hatg)$ corresponding
to~$\hatg$ is an associative algebra over~$\bc(q)$ with
generators~$E_i$, $F_i$, $K_i^{\pm1}$, $i\in\wh I$, $C^{\pm1/2}$ 
and~$D^{\pm1}$ subjects to the following relations
\begin{gather*}
\text{$C^{\pm1/2}$ are central and
$C=\prod_{i\in\wh I} K_i^{a_i}$, where~$\delta=
\displaystyle\sum_{i\in\wh I} a_i\alpha_i$}
\\
K_iK_i^{-1}=K_i^{-1}K_i=DD^{-1}=D^{-1}D=1,\quad K_iK_j=K_jK_i,\quad
\quad K_iD=DK_i,
\\
K_iE_j K_i^{-1}=q_i^{ a_{ij}}E_j,\qquad
K_iF_j K_i^{-1}=q_i^{-a_{ij}}F_j,\\
DE_j D^{-1}=q^{ \delta_{j,0}}E_j,\qquad
DF_j D^{-1}=q^{-\delta_{j,0}}F_j,\\
[E_i, F_j]=\delta_{i,j}\,\frac{K_i-K_i^{-1}}{q_i^{}-q_i^{-1}},
\intertext{}
\\
  \sum_{r=0}^{1-a_{ij}}(-1)^r
E_i^{(r)}E_{j}E_{i}^{(1-a_{ij}-r)}=0=
\sum_{r=0}^{1-a_{ij}}(-1)^r
F_{i}^{(r)}F_{j}F_{i}^{(1-a_{ij}-r)},
  \quad \text{if $i\ne j$}.
\end{gather*}
Let~$\bu_q^e$ be the quotient of~$\wh\bu_q$ by the two-sided ideal
generated by~$C^{\pm1/2}-1$. The algebra~$\bu_q$ is the subalgebra
of~$\bu_q^e$ generated by the~$E_i$, $F_i$ and~$K_i^{\pm1}$, $i\in\wh I$.

The elements~$E_i$, $F_i$ and~$K_i^{\pm1}$, $i\in I$
generate a subalgebra~$\bu_q^{fin}$ 
of~$\wh\bu_q$ which is isomorphic to the quantised 
enveloping algebra~$U_q(\lie g)$ of~$\lie g$. Notice also that,
for all~$i\in\wh I$ fixed, the elements~$E_i$, $F_i$ and~$K_i^{\pm1}$
generate a subalgebra of~$\wh\bu_q$ isomorphic to~$U_{q_i}(\lie{sl}_2)$.

\subl{P35}
One can introduce a $\bz$-grading on~$\wh\bu_q$ in the following way.
We say that~$x\in\wh\bu_q$ is homogeneous of degree~$k\in\bz$ 
if~$DxD^{-1}=q^k x$. That grading is obviously
well-defined since all generators of~$\wh\bu_q$ are homogeneous
and induces a $\bz$-grading on~$\bu_q$.
Given~$z\in\bc(q)^\times$, define an automorphism~$\phi_z$ of~$\wh\bu_q$
by~$\phi_z(x)=z^k x$ if~$x$ is homogeneous of degree~$k$. Evidently,
$\phi_z$ descends to an automorphism of~$\bu_q$.

Let~$M$ be a~$\wh\bu_q$ or~$\bu_q$-module. Denote by~$\phi_z^*M$ the
vector space~$M$ with the action of~$\wh\bu_q$ twisted by 
the automorphism~$\phi_z$, that is~$x \phi_z^*(m):=\phi_z(x)m$
for all~$x\in\bu_q$ or~$\wh\bu_q$, $m\in M$. Notice that the
map~$M\to\phi_z^*M$ is trivial as a map of vector spaces 
or~$\bu_q^{fin}$-modules.

Let~$M$ be a $\bu_q$-module. One can endow the loop space
$\wh M:=M\tensor_{\bc(q)} \bc(q)[t,t^{-1}]$ of~$M$ with the
structure of a $\wh\bu_q$-module by setting
$$
x(m\tensor f(t))=xm \tensor t^k f(t),\qquad
D^{\pm1}(m\tensor f(t))=m\tensor f(q^{\pm1}t),\qquad C^{\pm1/2}m=m,
$$
for all~$m\in M$, $f\in \bc(q)[t^{\pm1}]$ and for all~$x\in
\bu_q$ homogeneous of degree~$k$.

\subl{P36}
Let~$M$ be a $\bu_q$ (respectively, $\wh\bu_q$) module. We say that~$M$ is
a module of type~$1$ if~$M=\bigoplus_{\nu\in P_0} M_\nu$
(respectively, $M=\bigoplus_{\nu\in\wh P} M_\nu$), where
$M_\nu=\{ m\in M\,:\, K_i m=q_i^{\<\alpha_i^\vee,\nu>} m,\,
\forall\,i\in\wh I\}$ 
(respectively, $M_\nu=\{m\in M\,:\, K_i m=q_i^{\<\alpha_i^\vee,\nu>} m,\,
\forall\,i\in\wh I,\, Dm=q^{\<\partial,\nu>}m\}$). 
The subspaces~$M_\nu$ are called weight
subspaces of~$M$ and we call~$M$ {\em admissible} if~$\dim M_\nu<
\infty$ for all~$\nu\in P_0$ (respectively, for all~$\nu\in\wh P$). 
An element~$\nu\in P_0$ or~$\wh P$ is a weight of~$M$ if~$M_\nu\not=0$.

A module of type~$1$ is said to be of level~$k\in\bz$ if~$C$
acts on~$M$ by~$q^k\id$ and is said to be {\em integrable}
if the generators~$E_i$, $F_i$, $i\in\wh I$ act locally nilpotently
on~$M$. In other words, $M$ is a direct sum (possibly infinite),
of finite dimensional simple~$U_{q_i}(\lie{sl}_2)$-modules for all $i\in\wh I$.
Evidently, if~$M$ is a finite dimensional~$\bu_q$-module,
then~$\wh M$ is an integrable~$\wh\bu_q$-module. Moreover,
observe that all weights
of~$\wh M$ are of the form~$\nu+r\delta$ where~$\nu\in P_0$ and $r\in\bz$,
and that~$\wh M_{\nu+r\delta}$ is spanned by~$m\tensor t^r$ where~$m\in
M_{\nu}$. Thus, $\wh M$ is admissible.

\subl{P40}
It is well-known that~$\wh\bu_q$ admits a structure of 
a Hopf algebra. Throughout the
rest of this paper we will use the co-multiplication given
on generators by the following formulae
\begin{equation}\lbl{10}
\Delta(E_i)=E_i\tensor K_i^{-1}+1\tensor E_i,\quad
\Delta(F_i)=F_i\tensor 1+K_i\tensor F_i,
\end{equation}
the elements~$K_i^{\pm1}$, $D^{\pm1}$ and~$C^{\pm1/2}$ being group-like.
Then one can easily prove by induction on~$r$ that
\begin{equation}\lbl{20}
\begin{split}
&\Delta(E_i^{(r)})=\sum_{s=0}^r q_i^{-s(r-s)} E_i^{(s)}\tensor 
E_i^{(r-s)} K_i^{-s},\\
&\Delta(F_i^{(r)})=\sum_{s=0}^r q_i^{-s(r-s)} F_i^{(r-s)}K_i^s\tensor 
F_i^{(s)}
\end{split}
\end{equation}
Evidently, the above Hopf algebra structure descends to the algebra~$\bu_q$.
Henceforth, unless specified otherwise, a tensor product of two~$\wh\bu_q$
or~$\bu_q$ modules is assumed to be endowed with a structure of
a~$\wh\bu_q$ or~$\bu_q$ module with respect to the co-product~\loceqref{10}.

\subl{P45}
The algebras~$\wh\bu_q$ and~$\bu_q$ admit another presentation, known
as the Drinfel'd or loop-like presentation (cf.~\cite{Be,Dr,J}). Namely,
$\wh\bu_q$ is isomorphic to an associative algebra over~$\bc(q)$
generated by the~$x_{i,k}^\pm$, $h_{i,r}$, $K_i^{\pm1}$, $i\in I$,
$k\in\bz$, $r\in\bz\setminus\{0\}$, $C^{\pm1/2}$ and~$D^{\pm1}$
subjects to certain relations (see, for example, \cite{Be,C}).
Let us only mention that the~$x_{i,k}^{\pm}$ and the~$h_{i,k}$ are
homogeneous of degree~$k$ and~$x_{i,0}^{+}$ (respectively,
$x_{i,0}^-$) identifies with~$E_i$ (respectively, $F_i$).

For all $i\in I$ and~$r\in\bz$, 
define~$P_{i,\pm r}$ by equating the powers of~$u$
in the formal power series
$$
\sum_{r\ge0} P_{i,\pm r} u^r =  
\exp\Big(-\sum_{k>0} \frac{q^{\pm k}h_{i,\pm
k}}{[k]_i}u^k\Big).
$$
Then the~$P_{i,r}$, $i\in I$, $r\in\bz$ are homogeneous of degree~$r$
and generate the same subalgebra of~$\wh\bu_q$ as the~$h_{i,r}$.

A $\bu_q$-module $M$ is called $l$-highest weight with
highest weight~$(\lambda,\bpi^\pm)$, where~$\lambda\in P_0$ 
and~$\bpi^\pm=(\pi^\pm_1(u),\dots,\pi^\pm_\ell(u))$ with~$\pi^\pm_i(u)=
\sum_{k\ge 0} \pi_{i,\pm k} u^k\in\bc(q)[[u]]$ and $\pi_{i,0}=1$, 
if there exists a non-zero~$m\in M_\lambda$ such that~$M=\bu_q m$ 
and
$$
x_{i,k}^\pm m=0,\qquad P_{i,\pm k} m=\pi_{i,\pm k} m,
\qquad \forall i\in I,\, k\in\bz.
$$
Such an~$m$ is called an $l$-highest vector.
By~\cite{CP,CPweyl}, an $l$-highest weight module~$M$
with highest weight~$(\lambda,\bpi^\pm)$ is simple and finite dimensional
provided that $\bpi(u)=\bpi^+(u)=(\pi_1,\dots,\pi_\ell)$ 
is an $\ell$-tuple of polynomials,
$\deg \pi_i=\<\alpha_i^\vee,\lambda>$ and~$\pi_i^-(u)=
u^{\deg \pi_i} \pi_i(u^{-1})/(u^{\deg \pi_i} \pi_i(u^{-1})|_{u=0})$.
Moreover, all finite dimensional simple $\bu_q$-modules are
obtained that way. Henceforth we denote the simple
finite dimensional $l$-highest weight
module corresponding to an $\ell$-tuple~$\bpi$ of polynomials with
constant term~$1$ by~$V(\bpi)$. Let~$v_{\bpi}$ be the unique, up to
a scalar, $l$-highest weight vector of~$V(\bpi)$.

Let~$z\in\bc^\times$. Since the~$P_{i,k}$ are homogeneous of degree~$k$, 
$P_{i,\pm k} \phi_z^*(v_{\bpi})=z^{\pm k} \pi_{i,\pm k} v_{\bpi}$.
It follows that~$\phi_z^* V(\bpi)$ is isomorphic to~$V(\bpi_z)$
where~$\bpi_z(u)=\bpi(zu)$.

\subl{P50}
Let us conclude this section with a brief review of some facts about
crystals which we will need later. Throughout the rest of this paper, a
crystal is a set~$B$ endowed with maps~$e_i,f_i:B\to B\sqcup \{0\}$,
$\varepsilon_i,\varphi_i:B\to\bz$ for all~$i\in\wh I$ 
and~$\wt:B\to P$ or~$\wt:B\to\wh P$ 
satisfying the standard axioms
(see~\cite[1.2]{Ka1a} or~\cite[5.2]{JB}).
In particular, $\varphi_i(b)=\varepsilon_i(b)+\<\alpha_i^\vee,\wt b>$
for all~$b\in B$ and~$e_i, f_i$ for~$i$ fixed are quasi-inverses
of each other i.e. for all~$b,b'\in B$, $e_i b=b'$ if and only
if~$f_i b'=b$.
All crystals we consider are normal, that is~$\varepsilon_i(b)=
\max\{n\,:\, e_i^n b\in B\}$, $\varphi_i(b)=\max\{n\,:\, f_i^n b\in B\}$.
All morphisms of crystals will be assumed to be strict, that is,
commuting with all crystal operators. We say that~$B_1$ is a
subcrystal of~$B_2$ if there exists an injective morphism of
crystals~$B_1\to B_2$. If~$B_1$ is a subset of~$B_2$, we say that
$B_1$ is a subcrystal if the trivial embedding is a morphism of
crystals that is, if $B_1$ is a crystal with respect to the crystal
operations on~$B_2$ restricted to~$B_1$.

Let~$\Cm$ be the associative monoid generated by the operators~$e_i,
f_i\,:\,i\in\wh I$. A crystal~$B$ is generated by~$b\in B$ over~$\Cm$ 
if~$B=\Cm b:=\{x b\,:\, x\in\Cm\}\setminus
\{0\}$.
We say that a crystal~$B$ is indecomposable if it does not admit a 
non-empty subcrystal different from itself. By say~\cite[2.5]{G}
a crystal~$B$ is indecomposable if and only if $B$ is generated
by some~$b\in B$ over~$\Cm$. Moreover, if~$B=\Cm b$ for some~$b\in B$
then~$B=\Cm b'$ for all~$b'\in B$.

Given a family of crystals~$B_1,\dots,B_n$ one can introduce
a structure of a crystal on the set~$B_1\times\cdots\times B_n$,
which is called the tensor product of crystals
and denoted by~$B_1\tensor\cdots\tensor B_n$, in the following
way (cf.~\cite[1.3]{Ka1a}).
Given~$b=b_1\tensor\cdots \tensor b_n$, $b_i\in B_i$,
define the Kashiwara functions~$b\mapsto r_k^i(b)\,:\, i\in I$, $k\in\{1,
\dots, n\}$ by
$$
r_k^i(b)=\varepsilon_i(b_k)-\sum_{1\le j<k} \<\alpha_i^\vee,\wt b_j>.
$$
Then~$\varepsilon_i(b)$ is defined to be the maximal value of~$r_k^i(b)$
as a function of~$k$, $\wt b=\wt b_1+\cdots+\wt b_n$ and $e_i$
(respectively, $f_i$) acts in the leftmost (respectively, rightmost)
place where the maximal value of~$r_k^i(b)$ is attained. That
is known as Kashiwara's tensor product rule. It takes a particularly
nice form for~$n=2$ (cf.~\ref{ZCB20}).

Let~$B$ be a crystal with~$\wt:B\to P$. Its {\em affinisation} 
$\wh B=B\times \bz$ is a crystal with respect to the following operators.
Denote the pair~$(b,n)\in B\times\bz$ as~$b\tensor t^n$. Then
$\varepsilon_i(b\tensor t^n)=\varepsilon_i(b)$ and $\wt b\tensor t^n=
\wt b+n\delta\in\wh P$. Furthermore, if~$e_i b=0$, set~$e_i(b\tensor
t^n)=0$. Otherwise, $e_i(b\tensor t^n)=e_i b\tensor t^{n+\delta_{i,0}}$.
Similarly, if~$f_i b=0$, set~$f_i(b\tensor t^n)=0$. Otherwise,
set~$f_i(b\tensor t^n)=f_i b\tensor t^{n-\delta_{i,0}}$. This should
be regarded as the crystal analogue of the passage from a
$\bu_q$-module~$V$ to a~$\wh\bu_q$-module~$\wh V$ (cf.~\ref{P35}).

\section{General properties of $z$-crystal bases}

\subl{ZCB10}
Let~$M$ be an integrable $\bu_q$ or $\wh\bu_q$-module of type~$1$.
Fix~$i\in\wh I$ and let~$u$ be a weight vector of~$M$ of
weight~$\nu$. Then~$u$ can be written uniquely as
\begin{equation}\lbl{10}
u=\sum_{s\ge\max\{0,-\<\alpha_i^\vee,\nu>\}} F_i^{(s)} u_s,
\end{equation}
where~$u_s\in\ker E_i\cap M_{\nu+s\alpha_i}$ and~$u_s=0$ for~$s\gg0$.
The crystal operators of Kashiwara are defined as
\begin{equation}\lbl{20}
\tilde e_i u=\sum_{s\ge\max\{1,-\<\alpha_i^\vee,\nu>\}} F_i^{(s-1)} u_s,
\qquad
\tilde f_i u=\sum_{s\ge\max\{0,-\<\alpha_i^\vee,\nu>\}} F_i^{(s+1)} u_s.
\end{equation}
Observe that, since~$M$ is integrable, the operators~$\tilde e_i$,
$\tilde f_i$ are locally nilpotent.
\begin{defn}[{cf.~\cite[4.8]{CG}}]
Let~$z\in\bc^\times$. A $z$-crystal basis of~$M$ is a pair~$(L,B)$,
where~$L$ is a free~$\Ca$-submodule of~$M$ and~$B$ is a basis 
of~$\bc$-vector space $L/qL$ such that
\begin{enumerit}
\item $M=L\tensor_{\Ca} \bc(q)$.
\item $L=\bigoplus_{\lambda} L_\lambda$, $B=\coprod_{\lambda} B_\lambda$,
where~$L_\lambda=L\cap M_\lambda$ and 
$B_\lambda=B\cap (L_\lambda/q L_\lambda)$.
\item $L$ is preserved by the operators~$\tilde e_i$, $\tilde f_i$
for all~$i\in\wh I$. In particular, $\tilde e_i$, $\tilde f_i$ act
on~$L/qL$.
\item $\tilde e_i B,\tilde f_i B\subset z^{\bz\delta_{i,0}} B\cup\{0\}$,
for all~$i\in\wh I$.
\item For all $b,b'\in B$, $i\in\wh I$, 
$\tilde e_i b=z^{r\delta_{i,0}} b'$ if and only
if~$\tilde f_i b'=z^{-r\delta_{i,0}} b$.
\end{enumerit}
\end{defn}
For~$z=1$ the above definition reduces to Kashiwara's definition
of crystal bases (cf. for example~\cite{Ka1}).

Let~$(L,B)$ be a $z$-crystal basis of an integrable $\bu_q$ or
$\wh\bu_q$-module~$M$. Given~$b\in B$, let~$\varepsilon_i(b)$
(respectively, $\varphi_i(b)$) be the minimal non-negative integer~$n$
such that~$\tilde e_i^{n+1}b$ (respectively, $\tilde f_i^{n+1}b$)
equals zero modulo~$qL$. These are well-defined since $\tilde e_i$,
$\tilde f_i$ are locally nilpotent.
Furthermore, if~$b\in B_\lambda$, set~$\wt b=\lambda$.

We will need the following simple modification of~\cite[Lemma~20.1.2]{L}
(cf.~\cite[Lemma~4.8]{CG}).
\begin{lem}
Fix~$i\in\wh I$. Let~$u\in L_\lambda$ and 
write~$u=\sum_{s\ge\max\{0,-\<\alpha_i^\vee,\nu>\}} 
F_i^{(s)}u_s$ as in~\loceqref{10}.
Then
\begin{enumerit}
\item $F^{(r)}u_s\in L$ for all~$r,s\ge0$.
\item If $u\pmod{qL}\in B$ then there exists~$s_0$ such that
$u_s\in qL$, $s\not=s_0$, $u_{s_0}\pmod{qL}\in z^{\bz\delta_{i,0}} B$
and~$u=F^{(s_0)}u_{s_0}$.
\end{enumerit}
\end{lem}
\begin{pf}
The proof is an obvious modification of that of Lemma~20.1.2 in~\cite{L}.
\end{pf}

\subl{ZCB12}
Let~$(L,B)$ be a $z$-crystal basis of an integrable~$\bu_q$ 
or~$\wh\bu_q$-module~$M$. 
It follows immediately from~\defref{ZCB10} 
that the set~$\widetilde B=\bigcup_{r\in\bz} z^r B$ is a normal crystal.
Define an equivalence relation on~$\widetilde B$ by setting~$b\sim_z
b'$ if and only if~$b=z^r b'$ for some~$r\in\bz$. Then~$\widetilde B/
\!\sim_z$ identifies with~$B$ as a set. 
\begin{lem}
The set~$\widetilde B/\!\sim_z$ is a normal crystal with respect
to the operators~$\tilde e_i$, $\tilde f_i$, $\varepsilon_i$, $\varphi_i$, 
$i\in\wh I$ and $\wt$.
\end{lem}
\begin{pf}
Immediate.
\end{pf}
We call~$\widetilde B/\!\sim_z$ the crystal associated with~$B$.

\subl{ZCB15}
The following proposition justifies the definition
of $z$-crystal bases.
\begin{prop}
Let~$M$ be a finite-dimensional $\bu_q$-module and assume that~$M'=\phi_z^*M$
is not isomorphic to~$M$. Suppose that~$M$ admits a crystal basis
$(L,B)$. Then~$(L',B')$ where~$L'=\phi_z^*L$, $B'=\phi_z^*B\subset
L'/qL'$ is a $z$-crystal basis of~$M'$. Moreover, if~$b_1,b_2\in B$ such
that~$\tilde e_0 b_1 = b_2$ and~$b_1',b_2'$ are their images in~$B'$ then
$\tilde e_0 b_1'=z b'_2$.
\end{prop}
\begin{pf}
Since~$M\cong M'$ as~$\bu_q^{fin}$-module,
it is sufficient to verify~(iii)-(v) for~$i=0$. 

Let~$u\in L_\nu$ be a weight vector and let~$u_s$, $s\ge0$
be as in~\eqref{ZCB10.10}. Set~$u'=\phi_z^* u$. Then
$$
u'=\sum_{s\ge\max\{0,-\<\alpha_0^\vee,\nu>\}} z^{s} F_0^{(s)} \phi_z^* u_s=
\sum_{s\ge\max\{0,-\<\alpha_0^\vee,\nu>\}} F_0^{(s)} u'_s,
$$
where~$u'_s=z^{s} \phi_z^* u_s$. This provides the unique 
decomposition of the form~\eqref{ZCB10.10} for~$u'\in L'$. Then
$$
\phi_z^*(\tilde e_0 u)=\sum_{s\ge\max\{1,-\<\alpha_0^\vee,\nu>\}} 
\mskip-25mu z^{s-1} F_0^{(s-1)} 
\phi_z^*(u_s)=z^{-1}\mskip-10mu\sum_{s\ge\max\{1,-\<\alpha_0^\vee,\nu>\}} 
\mskip-25mu
F_0^{(s-1)} u'_s=z^{-1}\tilde e_0 u'.
$$
Since~$\tilde e_0 u\in L$, it follows that~$\tilde e_0 u'=z
\phi_z^*(\tilde e_0 u)\in L'$. Similarly, $\tilde f_0 u'=
z^{-1} \phi_z^*(\tilde f_0 u)\in L'$. Since $u_s\in L$ for all~$s\ge0$
by~\lemref{ZCB10}(i), it follows that~$u'_s\in L'$ for 
all~$s\ge0$.

Furthermore, suppose that~$b=u\pmod{qL}\in B$. Then
by~\lemref{ZCB10}(ii) there exists~$s_0$ such that~$u_s\in qL$,
$s\not=s_0$, $u_{s_0}\pmod{qL}\in z^{\bz} B$ and~$b=F_0^{(s_0)}u_{s_0}
\pmod{qL}$. Moreover, $\tilde e_0 b=F_0^{(s_0-1)}u_{s_0}\pmod{qL}$,
$\tilde f_0 b=F_0^{(s_0+1)} u_{s_0}\pmod{qL}$. Let~$b'=\phi_z^*b=
\phi_z^*u\pmod{qL'}\in B'$. It follows immediately that~$\tilde e_0b'=
z\phi_z^*(\tilde e_0 u)\pmod{qL'}\in z^\bz B'\cup\{0\}$,
$\tilde f_0b'=z^{-1} \phi_z^*(\tilde f_0 u)\pmod{qL'}\in z^\bz B'\cup\{0\}$.

Finally, suppose that~$\tilde e_0 b=b_1\in B$ and let~$b_1'=\phi_z^* b_1$.
Then, as above, $\tilde e_0 b'=z b_1'$. On the other hand, $b_1'=
z^{s_0-1} F^{(s_0-1)} \phi_z^*(u_{s_0})\pmod{qL'}$, whence
$\tilde f_0 b_1'=z^{s_0-1} F^{(s_0)}\phi_z^*(u_{s_0})=z^{-1}b'$.
\end{pf}
\begin{rem}
Similarly, one can prove that if $(L,B)$ is a $z$-crystal basis of~$M$
and~$\phi_z^*M$ is not isomorphic to~$M$ 
then $(\phi_z^*L,\phi_z^*B)$ is a $z$-crystal basis of $\phi^*_z M$.
Moreover, if~$\tilde e_0 b=z^k b_1$ for some~$b,b_1\in B$ and
$b'$, $b_1'$ denote their images in~$\phi_z^*B$ then 
$\tilde e_0 b'=z^{k+1} b_1'$.
\end{rem}

\subl{ZCB17}
The following Lemma is rather standard (cf.~\cite[Corollary~17.4.2]{L}).
We deem it necessary to present its proof here since the argument 
in~\cite{L} is based on the use of Kashiwara's bilinear form and cannot be
modified for $z$-crystal bases.
\begin{lem}
Let~$M_i$, $i=1,2$ be finite dimensional $\bu_q(\lie{sl}_2)$~modules 
of type~$1$ and 
fix~$v_i$, $i=1,2$ such that~$Kv_i=q^{t_i}v_i$, $t_i\ge0$ and~$Ev_i=0$.
Let~$\Cl_i$ be the~$\Ca$-module generated by the $F^{(s)}v_i$, $0\le s\le t_i$.
Then
\begin{enumerit}
\item
$\Cl:=\Cl_1\tensor_{\Ca}\Cl_2$ is preserved by the operators~$\tilde e$,
$\tilde f$ acting on the module~$M_1\tensor M_2$.
\item
There exist unique, up to multiplication by an element
of~$1+q\Ca$, $u_r\in\ker E\cap\Cl$, $0\le r\le\min\{t_1,t_2\}$ such
that $u_r\notin q\Cl$, $Ku_r=q^{t_1+t_2-2r}u_r$, $\Cl$ is a direct sum
of $\Ca$-modules generated by~$F^{(b)}u_r$, $0\le b\le t_1+t_2-r+1$ and
for all~$0\le s_i\le t_i$, $i=1,2$, there exists a unique~$s$,
$0\le s\le \min\{t_1,t_2\}$ such that~$F^{(s_1)}v_1\tensor F^{(s_2)}v_2=
F^{(s_1+s_2-s)}u_s\pmod{q\Cl}$.
\item For all~$0\le s_i\le t_i$,
\begin{alignat*}{2}
&\tilde e(F^{(s_1)}v_1\tensor F^{(s_2)}v_2)=F^{(s_1-1)}v_1\tensor 
F^{(s_2)}v_2\pmod{q\Cl},&\qquad& t_1\ge s_1+s_2\\
&\tilde e(F^{(s_1)}v_1\tensor F^{(s_2)}v_2)=F^{(s_1)}v_1\tensor 
F^{(s_2-1)}v_2\pmod{q\Cl},&& t_1< s_1+s_2\\
&\tilde f(F^{(s_1)}v_1\tensor F^{(s_2)}v_2)=F^{(s_1+1)}v_1\tensor 
F^{(s_2)}v_2\pmod{q\Cl},&\qquad& t_1> s_1+s_2\\
&\tilde f(F^{(s_1)}v_1\tensor F^{(s_2)}v_2)=F^{(s_1)}v_1\tensor 
F^{(s_2+1)}v_2\pmod{q\Cl},&\qquad& t_1\le s_1+s_2,
\end{alignat*}
where~$F^{(s)}v_i=0$ if~$s<0$.
\end{enumerit}
\end{lem}
\begin{pf}
Let~$V(n)$, $n\ge0$ denote the unique~$(n+1)$-dimensional 
simple~$\bu_q(\lie{sl}_2)$ module.
It is sufficient to prove the Lemma for~$M_i\cong V(t_i)$.
The argument is by induction on~$t_1$ and is rather standard.

$1^\circ\mskip-7mu.$~Suppose first 
that~$t_1=0$. Then~$F^{(s)}(v_1\tensor v_2)=v_1\tensor F^{(s)}
v_2$ and~$E(v_1\tensor v_2)=0$. The proposition is then trivial.

$2^\circ\mskip-7mu.$~Suppose that~$t_1=1$ and set~$u_0=v_1\tensor v_1$,
$$
u_1=v_1\tensor Fv_2-q^{t_2}[t_2]_q Fv_1\tensor v_2=
v_1\tensor Fv_2-q\,\frac{1-q^{2t_2}}{1-q^{2\hphantom{t_2}}}\,Fv_1\tensor v_2.
$$
Then~$u_0$, $u_1$ generate~$\ker E\cap\Cl$ and~$u_1=
v_1\tensor Fv_2\pmod{q\Cl}$.  Furthermore,
\begin{align}
&F^{(b)}u_0=q^b v_1\tensor F^{(b)}v_2+Fv_1\tensor F^{(b-1)}v_2\lbl{10}\\
&F^{(b)}u_1=\frac{1-q^{2(b+1)}}{1-q^{2\hphantom{(b+1)}}}\,
v_1\tensor F^{(b+1)}v_2+q\,\frac{q^{2t_2-b}-q^{b}}{1-q^2}
Fv_1\tensor F^{(b)}v_2,\lbl{20}
\end{align}
where we used~\eqref{P40.20}.
Since~$F^{(b)}v_2=0$ if~$b>t_2$, it follows immediately that
$F^{(b)}u_0,F^{(b)}u_1\in\Cl$ for all~$b\ge0$. Moreover, by
the above formulae, $v_1\tensor F^{(b)}v_2=F^{(b-1)}u_1\pmod{q\Cl}$
whilst~$Fv_1\tensor F^{(b-1)}v_2=F^{(b)}u_0\pmod{q\Cl}$, $b>0$. 
It follows that the matrix of~$F^{(b)}u_0$, $F^{(b-1)}u_1$ in
the basis of~$Fv_1\tensor F^{(b-1)}v_2$, $v_1\tensor F^{(b)}v_2$ is
diagonal and its diagonal entries equal~$1\pmod{q\Ca}$. Therefore,
that matrix is invertible over~$\Ca$ and so the $F^{(b)}u_0$,
$F^{(b-1)}u_1$ and $Fv_1\tensor F^{(b-1)}v_2$, $v_1\tensor F^{(b)}v_2$
generate the same~$\Ca$-module which completes the proof of~(ii).

Since~$Eu_0=0$, it follows that~$\tilde e(v_1\tensor v_2)=0$. On the other
hand, $\tilde f(v_1\tensor v_2)=F u_0=Fv_1\tensor v_2\pmod{q\Cl}$, which
agrees with the formulae in~(iii).

Suppose now that~$b>0$. By the above,
$F^{(s)}v_1\tensor F^{(b-s)}v_2=x_s F^{(b)}u_0+y_s F^{(b-1)}u_1$,
where~$x_s\in\delta_{s,1}+q\Ca$, $y_s\in\delta_{s,0}+q\Ca$. Then, by 
definition of Kashiwara's operators,
$$
\tilde e(F^{(s)}v_1\tensor F^{(b-s)}v_2)=x_s F^{(b-1)}u_0+
y_s F^{(b-2)}u_1,
$$
In particular, $\tilde e$ preserves~$\Cl$.
If~$s=0$ then the above expression equals~$F^{(b-2)}u_1\pmod{q\Cl}=
v_1\tensor F^{(b-1)}v_2\pmod{q\Cl}$ provided that~$b\ge2$ (and so
$s_1+s_2=b>t_1$), which agrees with the formulae in~(iii).
If~$b=1$ (that is,~$s_1+s_2\le t_1$), 
$\tilde e(v_1\tensor Fv_2)=x_s u_0=0\pmod{q\Cl}$ as expected. 
Similarly, if~$s=1$, we get
$$
\tilde e(Fv_1\tensor F^{(b-1)}v_2)= F^{(b-1)}u_0\pmod{q\Cl}=Fv_1\tensor
F^{(b-2)}v_2\pmod{q\Cl},
$$
as desired. The formulae for the action of~$\tilde f$ are proved similarly.

$3^\circ\mskip-7mu.$~Suppose 
that~(i)--(iii) are proved for all~$t_1\le t$, $t>0$. 
It is well-known (cf., for example, \cite[4.3]{JB}) 
that~$V(t+1)$ can be realised as a simple
submodule of~$V(1)\tensor V(t)$ generated by the tensor product of the
corresponding highest weight vectors.  Thus, we can write~$v_1$ from
the assertion of the Lemma as~$v''_1\tensor v'_1$, where~$Ev_1'=Ev_1''=0$,
$F^2v''_1=0$, $F^{t+1}v_1'=0$. Let~$\Cl'$ be the $\Ca$-module generated
by the~$F^{(s_1)}v'_1\tensor F^{(s_2)}v_2$ and denote by~$u'_r$,
$0\le r\le\min\{t,t_2\}$ the 
elements of~$\ker E\cap \Cl'$ satisfying~$Ku'_r=q^{t+t_2-2r}u'_r$ given
by the induction hypothesis. Let~$\Cl''=\Ca v_1''+\Ca Fv_1''$. 
It follows from~\loceqref{20} that
the $\Ca$-module~$\Cl$ generated by the~$F^{(s_1)}v_1\tensor F^{(s_2)}v_2$ is
an $\Ca$-submodule of~$\Cl''\tensor_{\Ca}\Cl'$. 

By the induction hypothesis, $\Cl'=\bigoplus_r 
\Big(\bigoplus_b \Ca F^{(b)}u'_r\Big)$.
Applying the second part of the proof to~$\Cl''\tensor_{\Ca} \bigoplus_b
\Ca F^{(b)}u'_r$, $0\le r\le\min\{t,t_2\}$, 
we conclude that~$\tilde e$, $\tilde
f$ preserve $\Cl''\tensor_{\Ca}\Cl'$. Since~$\Cl$ is contained in the
intersection of~$\Cl''\tensor_\Ca\Cl'$ with a submodule of~$V(1)\tensor
V(t)\tensor V(t_2)$, it follows that $\tilde e$, $\tilde f$
preserve $\Cl$.

The next step is to prove the formulae in~(iii). Since~$\tilde e$, $\tilde f$
preserve~$\Cl$, we can do all the computations modulo~$q\Cl$.

Consider first~$v_1\tensor
F^{(s_2)}v_2=v''_1\tensor v'_1\tensor F^{(s_2)}v_2$, $s_2\le t_2$. 
By the induction hypothesis, $v'_1\tensor F^{(s_2)}v_2=
F^{(s_2-s)}u'_s\pmod{q\Cl'}$ for some~$0\le s\le \min\{t,t_2\}$.
Suppose first that~$s_2\le t$. Then~$\tilde e(v'_1\tensor F^{(s_2)}v_2)=0
\pmod{q\Cl'}=F^{(s_2-s-1)}u'_s$
by the induction hypothesis, whence~$s_2=s$. It follows that
$\tilde e(v_1\tensor F^{(s_2)}v_2)=\tilde e(v''_1\tensor u'_{s_2})
\pmod{q\Cl}=0$, as desired. Suppose that~$s_2=t+k$, $k>0$. 
Then~$v_1'\tensor F^{(s_2)}v_2=\tilde f^k(v_1'\tensor F^{(t)}v_2)
\pmod{q\Cl'}=F^{(k)}u'_t\pmod{q\Cl'}$ by the induction hypothesis.
Then~$\tilde e(v_1\tensor F^{(s_2)}v_2)=\tilde e(v''_1\tensor 
F^{(k)}u'_t)\pmod{q\Cl}$. By the first part of the proof, the
latter expression equals zero if~$k=1$ (that is, $s_2=t+1$) and
$v''_1\tensor F^{(k-1)}u'_t=v_1\tensor F^{(s_2-1)}v_2\pmod{q\Cl}$
if~$k>1$ (that is, $s_2>t+1$). Both agree with the formulae in~(iii). 

It follows from~\loceqref{10} that~$F^{(s_1)}v_1=Fv''_1\tensor
F^{(s_1-1)}v'_1\pmod{q(\Cl''\tensor_{\Ca}\Cl')}$, $0<s_1\le t+1$. 
Suppose first that~$s_2\le t$. Then
$\tilde f(v''_1\tensor u'_{s_2})=Fv''_1\tensor u'_{s_2}\pmod{q\Cl}
=Fv''_1\tensor v'_1\tensor F^{(s_2)}v_2=Fv_1\tensor F^{(s_2)}v_2
\pmod{q\Cl}$ with agrees with~(iii). 
Similarly, if~$s_2>t$, $\tilde f(v_1\tensor
F^{(s_2)}v_2)=\tilde f(v''_1\tensor F^{(s_2-t)}u'_t)\pmod{q\Cl}=
v''_1\tensor F^{(s_2-t+1)}u'_t\pmod{q\Cl}$ by the second
part of the proof. Thus, $\tilde f(v_1\tensor F^{(s_2)}v_2)=
v_1\tensor F^{(s_2+1)}v_2\pmod{q\Cl}$ as desired.

Consider now~$F^{(s_1)}v_1\tensor F^{(s_2)}v_2$ with~$0<s_1\le t+1$.
Using the induction hypothesis, we get
\begin{alignat*}{2}
\tilde e(F^{(s_1)}v_1\tensor F^{(s_2)}v_2)&=\tilde e(Fv''_1\tensor
F^{(s_1-1)}v'_1\tensor F^{(s_2)}v_2)&&\pmod{q\Cl}\\&=
\tilde e(Fv''_1\tensor F^{(s_1+s_2-s-1)}u'_s)&&\pmod{q\Cl}
\end{alignat*}
for some~$s$, $0\le s\le\min\{t,t_2\}$.
Suppose first that~$s_1+s_2-s\le 1$. Then, by the second
part of the proof,
$$
\tilde e(Fv''_1\tensor F^{(s_1+s_2-s-1)}u'_s)=
v''_1\tensor F^{(s_1+s_2-s-1)}u'_s\pmod{q\Cl}.
$$
Yet~$s_1+s_2-s\ge 1$, hence~$F^{(s_1-1)}v'_1\tensor F^{(s_2)}v_2=u_{s_1+s_2-1}
\pmod{q\Cl'}$. In particular, $\tilde e(F^{(s_1-1)}v'_1\tensor
F^{(s_2)}v_2)=0\pmod{q\Cl'}$.
Suppose that~$s_1+s_2\le t+1$. Then, by the induction hypothesis
$0=\tilde e(F^{(s_1-1)}v'_1\tensor F^{(s_2)}v_2)=F^{(s_1-2)}v'_1\tensor
F^{(s_2)}v_2\pmod{q\Cl'}$. We conclude that~$s_1=1$. Thus
$$
\tilde e(F^{(s_1)}v_1\tensor F^{(s_2)}v_2)=v''_1\tensor v'_1\tensor 
F^{(s_2)}v_2\pmod{q\Cl}=v_1\tensor F^{(s_2)}v_2\pmod{q\Cl},
$$
which agrees with the formulae in~(iii). On the other hand, if~$s_1+s_2>t+1$,
then, by the induction hypothesis, $\tilde e(F^{(s_1-1)}v'_1\tensor
F^{(s_2)}v_2)=F^{(s_1-1)}v'_1\tensor F^{(s_2-1)}v_2$, whence~$s_2=0$ 
and~$s_1>t+1$ which is a contradiction.

Finally, assume that~$s_1+s_2-s>1$. Then, by the second part of the proof,
$\tilde e(Fv''_1\tensor F^{(s_1+s_2-s-1)}u'_s)=
Fv''_1\tensor F^{(s_1+s_2-s-2)}u'_s\pmod{q\Cl}$. Yet, by the induction
hypothesis, $F^{(s_1+s_2-s-2)}u'_s=\tilde e(F^{(s_1-1)}v'_1\tensor 
F^{(s_2)})\pmod{q\Cl'}$. The latter expression equals modulo
$q\Cl'$, by the induction
hypothesis, $F^{(s_1-2)}v'_1\tensor F^{(s_2)}v_2$ if~$s_1+s_2-1\le t$
and~$F^{(s_1-1)}v'_1\tensor F^{(s_2-1)}v_2$ otherwise. Thus,
\begin{alignat*}{2}
\tilde e(F^{(s_1)}v_1\tensor F^{(s_2)}v_2)&=Fv''_1\tensor F^{(s_1-2+k)}v'_1
\tensor F^{(s_2-k)}v_2&&\pmod{q\Cl}\\&=
F^{(s_1-1+k)}v_1\tensor F^{(s_2-k)}v_2&&\pmod{q\Cl},
\end{alignat*}
where~$k$ equals zero if~$s_1+s_2\le t+1$ and~$1$ otherwise. That proves
the first two formulae in~(iii). In order to prove the last two formulae,
observe that, since~$s_1+s_2-s\ge1$,
$\tilde f(Fv''_1\tensor F^{(s_1+s_2-s-1)}u'_s)=Fv''_1\tensor F^{(s_1+s_2-s)}
u'_s\pmod{q\Cl}=Fv''_1\tensor \tilde f(F^{(s_1-1)}v'_1\tensor
F^{(s_2)}v_2)\pmod{q\Cl}$. It remains to apply the induction hypothesis.

The last step is to prove~(ii). Set, for~$0\le r\le\min\{t+1,t_2\}$,
$$
u_r = \sum_{a=0}^r c_{r,a} F^{(a)}v_1\tensor F^{(r-a)} v_2,
$$
where~$c_{r,0}=1$ and, for~$1\le a\le r$,
$$
c_{r,a} = (-1)^a \prod_{j=1}^a q^{t_2-2(r-j)}\,\frac{[t_2-r+j]_q}{[t-j+2]_q}
=(-1)^a q^{a(t-r+2)}\prod_{j=1}^a \frac{1-q^{2(t_2-r+j)}}{1-q^{2(t-j+2)}}.
$$
Then~$Eu_r=0$ and~$Ku_r=q^{t+t_2+1-2r}u_r$. Evidently, $u_r\in\Cl$ and
$u_r = v_1\tensor F^{(r)}v_2\pmod{q\Cl}$. We claim that, for all $0\le s_1\le
t+1$, $0\le s_2\le t_2$ there exist a unique $0\le s\le\min\{t+1,t_2\}$
such that~$F^{(s_1)}v_1\tensor F^{(s_2)}v_2=\tilde 
f^{s_1+s_2-s}u_s\pmod{q\Cl}=F^{(s_1+s_2-s)}u_s\pmod{q\Cl}$. Evidently, (ii)
follows immediately from the claim.

In order to prove the claim, observe first that~$v_1\tensor F^{(s_2)}v_2=
u_{s_2}\pmod{q\Cl}$, $0\le\min\{t+1,t_2\}$. If~$t_2\le t+1$ that gives
$v_1\tensor F^{(s_2)}v_2$ for all~$0\le s_2\le t_2$. Otherwise, by~(iii),
$\tilde f^k(v_1\tensor F^{(t+1)}v_2)=v_1\tensor F^{(t+k+1)}v_2\pmod{q\Cl}$.
Thus, $v_1\tensor F^{(s_2)}v_2=\tilde f^{s_2-t-1} u_{t+1}$, $s_2>t+1$.
Consider further~$F^{(s_1)}v_1\tensor F^{(s_2)}v_2$, $s_1>0$. We use induction
on~$s_1$. If~$s_1+s_2<t+1$
then we have, by~(iii), $F^{(s_1+1)}v_1\tensor F^{(s_2)}v_2=
\tilde f(F^{(s_1)}v_1\tensor F^{(s_2)}v_2\pmod{q\Cl}=\tilde f^{s_1+s_2-s+1}
u_s\pmod{q\Cl}$, where~$s$ is such that~$F^{(s_1)}v_1\tensor F^{(s_2)}v_2=
\tilde f^{s_1+s_2-s}u_s\pmod{q\Cl}$. Finally, suppose that
$s_1+s_2=t+1+k$, $k\ge0$. We may assume that~$s_1<t+1$ for otherwise
$F^{(s_1+1)}v_1=0$. Set~$l=t+1-s_1>0$. Then~$s_2=k+l\ge l$ 
and~$F^{(s_1)}v_1\tensor F^{(l-1)}v_2=
\tilde f^{s_1+l-s-1}u_s\pmod{q\Cl}$ by the induction hypothesis. Using~(iii)
repeatedly we conclude that~$F^{(s_1+1)}v_1\tensor F^{(s_2)}v_2=
\tilde f^{s_2-l+2}(F^{(s_1)}v_1\tensor F^{(l-1)}v_2)\pmod{q\Cl}
=\tilde f^{s_1+s_2-s+1}u_s\pmod{q\Cl}$, which completes the proof of 
the claim.
\end{pf}
\subl{ZCB20}
Let $M_i$, $i=1,2$ be finite dimensional $\bu_q$-modules or admissible
integrable $\wh{\bu}_q$-modules. Suppose 
that $M_1$ admits a crystal basis $(L_1,B_1)$ 
and that~$M_2$ admits a $z$-crystal
basis $(L_2,B_2)$ for some~$z\in\bc^\times$. 
\begin{prop}
The pair $(L,B)$, where~$L=L_1\tensor_{\Ca} L_2$ and~$B=\{b_1\tensor b_2\,:
b_i\in B_i\}$, is a $z$-crystal basis of~$M_1\tensor M_2$.
Moreover, for all~$b_i\in B_i$, $i=1,2$
\begin{align*}
&\tilde e_i(b_1\tensor b_2)=\begin{cases}
\tilde e_ib_1\tensor b_2,\quad&\varphi_i(b_1)\ge \varepsilon_i(b_2)\\
b_1\tensor \tilde e_ib_2,&\varphi_i(b_1)< \varepsilon_i(b_2)
\end{cases}
\\
&\tilde f_i(b_1\tensor b_2)=\begin{cases}
\tilde f_ib_1\tensor b_2,\quad&\varphi_i(b_1)> \varepsilon_i(b_2)\\
b_1\tensor \tilde f_ib_2,&\varphi_i(b_1)\le \varepsilon_i(b_2).
\end{cases}
\end{align*}
\end{prop}
\begin{pf}
The proof is essentially the same as that of~\cite[Theorem~20.2.2]{L}. 
We only have to verify the properties of a $z$-crystal basis for~$i=0$.
Set
$$
G^t_i:=\{ v \in L_i\,:\, K_0v=q_0^t v,\,
E_0v=0,\, v\in B_i\pmod{qL_i}\}.
$$
Then, by Lemma~\lemref{ZCB10}, 
$F_0^{(s_1)}v_1\in B_1\pmod{qL_1}$ for all~$v_1\in G^{t_1}_1$ 
and~$0\le s_1\le t_1$
and all elements of~$B_1$ are obtained that way. Similarly, 
for all~$v_2\in G^{t_2}_2$ and~$0\le s_2\le t_2$ there exists~$r=
r(v_2,s_2)\in\bz$
such that~$z^r F_0^{(s_2)}v_2\in B_2\pmod{qL_2}$ 
and all elements of~$B_2$ are obtained that way.
Since the weight spaces of $M_i$, $i=1,2$ are finite-dimensional, 
it follows by Nakayama's Lemma that
the~$\Ca$-module $L_i$ is generated over~$\Ca$ 
by the~$F_0^{(s_i)}v_i$, $v_i\in
G_i^{t_i}$, $0\le s_i\le t_i$.
Therefore, $L$ is generated over~$\Ca$ by the~$F_0^{(s_1)}v_1\tensor
F_0^{(s_2)}v_2$, $v_i\in G_i^{t_i}$, $0\le s_i\le t_i$, $i=1,2$.
Using~\lemref{ZCB17}(iii) we conclude that 
$\tilde e_0$, $\tilde f_0$ map the generators of the~$\Ca$-module~$L$ into~$L$
and hence act on~$L$. The rest of the properties of a
$z$-crystal basis and Kashiwara's tensor product
rule follows readily from~\lemref{ZCB17}(iii). 
\end{pf}

\subl{ZCB40}
Let~$V$ be a finite-dimensional simple $\bu_q$-module and assume that
$V$ admits a $z$-crystal basis~$(L,B)$ for some~$z\in\bc^\times$.
Let~$\wh V$ be as in~\ref{P35}
\begin{lem}
Set~$\wh L=L\tensor_{\Ca}\Ca[t,t^{-1}]$, $\wh B=\{b\tensor t^r\,:\,
b\in B,\,r\in\bz\}$. Then~$(\wh L,\wh B)$ is a $z$-crystal basis of the
$\wh\bu_q$-module~$\wh V$. Moreover, for all~$b\in B$, $r\in\bz$,
$$
\tilde e_i(b\tensor t^r)=(\tilde e_i b)\tensor t^{r+\delta_{i,0}}
\pmod{q\wh L},
\qquad \tilde f_i(b\tensor t^r)=(\tilde f_i b)\tensor t^{r-\delta_{i,0}}
\pmod{q\wh L}.
$$
In other words, the associated crystal of~$\wh B$ is the affinisation
of the associated crystal of~$B$ in the sense of the definition
given in~\ref{P50}.
\end{lem}
\begin{pf}
Take~$u\in L$ of weight~$\lambda$ and write 
$u=\sum_{s\ge\max\{0,-\<\alpha_i^\vee,\lambda>\}} F_i^{(s)} u_s$
as in~\eqref{ZCB10.10}. Evidently,
$$
u\tensor t^r=\sum_{s\ge\max\{0,-\<\alpha_i^\vee,\lambda>\}} F_i^{(s)} 
(u_s\tensor t^{r+s\delta_{i,0}}),
$$
which is the decomposition~\eqref{ZCB10.10} for~$u\tensor t^r$.
Then by the definition of Kashiwara's operators,
\begin{align*}
\tilde e_i(u\tensor t^r)&=\sum_{s\ge\max\{1,-\<\alpha_i^\vee,\lambda>\}}
F_i^{(s-1)}(u_s\tensor t^{r+s\delta_{i,0}})\\&=
\Big(\sum_{s\ge\max\{1,-\<\alpha_i^\vee,\lambda>\}}
F_i^{(s-1)} u_s\Big)\tensor t^{r+\delta_{i,0}}=(\tilde e_i u)\tensor t^{r
+\delta_{i,0}}\\
\tilde f_i(u\tensor t^r)&=\sum_{s\ge\max\{0,-\<\alpha_i^\vee,\lambda>\}}
F_i^{(s+1)}(u_s\tensor t^{r+s\delta_{i,0}})\\&=
\Big(\sum_{s\ge\max\{0,-\<\alpha_i^\vee,\lambda>\}}
F_i^{(s+1)} u_s\Big)\tensor t^{r-\delta_{i,0}}=(\tilde f_i u)\tensor t^{r
-\delta_{i,0}}.
\end{align*}
The assertion follows immediately from the above formulae and the properties
of a $z$-crystal basis.
\end{pf}

\section{Quantum loop modules and their $z$-crystal bases}

\subl{ZCB30}
Let~$\bpi^0$ be an $\ell$-tuple of polynomials over~$\bc(q)$
with constant term~$1$ and
suppose that~$\bpi^0(zu)\not=\bpi^0(u)$, $z\in\bc^\times$ 
as a set of polynomials. Given $\ell$-tuples of polynomials~$\bpi=(\pi_i)_{
i\in I}$,
$\bpi'=(\pi'_i)_{i\in I}$ set~$\bpi\bpi'=(\pi_i\pi_i')_{i\in I}$.

Retain the notations of~\ref{P45} and
suppose that the finite dimensional $\bu_q$-module $V(\bpi^0)$ admits
a crystal basis~$(L(\bpi^0),B(\bpi^0))$. Fix~$m\in\bn$ which does not
exceed the multiplicative order of~$z$ and set~$\bpi=
\bpi^0\bpi^0_z\cdots\bpi^0_{z^{m-1}}$. Then
$V(\bpi)$ is isomorphic to~$V(\bpi^0)\tensor V(\bpi^0_z)
\tensor\cdots\tensor V(\bpi^0_{z^{m-1}})$ by~\cite{C}.
Furthermore, set~$L(\bpi)=L(\bpi^0)\tensor_{\Ca}\phi^*_z L(\bpi^0)
\tensor_\Ca\cdots\tensor_\Ca \phi^*_{z^{m-1}}L(\bpi^0)$ and define
$B(\bpi)$ accordingly. Since~$\phi_z^*$ is the identity map on the level
of vector spaces, $B(\bpi)$ identifies with~$B(\bpi^0)^{\tensor m}=\{b_1\tensor
\cdots \tensor b_m\,:\,b_i\in B(\bpi^0)\}$. 
\begin{prop}
The pair~$(L(\bpi), B(\bpi))$ is a $z$-crystal basis of~$V(\bpi)$. Moreover,
for all~$b_1,\dots,b_m\in B(\bpi^0)$,
\begin{alignat*}{2}
&\tilde e_i (b_1\tensor\cdots\tensor b_m)=
z^{r-1}& &b_1\tensor\cdots\tensor b_{r-1}\tensor \tilde e_i b_r \tensor b_{r+1}
\tensor\cdots\tensor b_m\\
&\tilde f_i (b_1\tensor\cdots\tensor b_m)=
z^{-s+1}& &b_1\tensor\cdots\tensor b_{s-1}\tensor \tilde e_i b_s \tensor b_{s+1}
\tensor\cdots\tensor b_m,
\end{alignat*}
where~$r$ and~$s$ are determined by Kashiwara's tensor product rule. 
In particular, the associated crystal of~$B(\bpi)$ is isomorphic
to~$B(\bpi^0)^{\tensor m}$.
\end{prop}
\begin{pf}
The proof is by induction on~$m$, the induction base being trivial.
Recall that~$V(\bpi^0_a)=\phi^*_a V(\bpi^0)$. 
Set~$V_k=V(\bpi^0)\tensor V(\bpi^0_z)\tensor\cdots\tensor V(\bpi^0_{z^{k-1}})$,
$k>0$ and define~$L_k$, $B_k$ accordingly. Suppose that~$(L_k,B_k)$ is
a $z$-crystal basis for~$V_k$. Then~$V_{k+1}\cong V_1\tensor
\phi_z^*V_k$ and~$(\phi^*_z L_k,\phi^*_z B_k)$ is a $z$-crystal basis
of~$V_k$ by~\remref{ZCB15}. Then~$(L_1\tensor \phi^*_z L_k, B_1\tensor
\phi^*_z B_k)=(L_{k+1},B_{k+1})$ is a $z$-crystal basis of~$V_k$
by~\propref{ZCB20}. The formulae follow immediately from these 
in~\propref{ZCB20}.
\end{pf}

\subl{ZCB60}
Let~$\zeta$ be an $m$th primitive root of unity. Let~$\bpi^0$ be a tuple
of polynomials such that~$\bpi^0(\zeta u)\not=\bpi^0(u)$ as a set
of polynomials. Fix an $l$-highest weight vector~$v_{\bpi^0}$ in~$V(\bpi^0)$
and write~$v_{\bpi^0_z}=\phi_z^*v_{\bpi^0}$.
Let~$V(\bpi)=V(\bpi^0)\tensor V(\bpi^0_\zeta)\tensor \cdots\tensor 
V(\bpi^0_{\zeta^{m-1}})$ and set~$v_{\bpi}=
v_{\bpi^0}\tensor v_{\bpi^0_\zeta}\tensor \cdots\tensor
v_{\bpi^0_{\zeta^{m-1}}}$.
By~\cite{C}, $V(\bpi)$ is a simple $\bu_q$-module
and there exists a unique isomorphism of~$\bu_q$-modules
$$
\tau:V(\bpi)\to V(\bpi^0_{\zeta^{m-1}})\tensor V(\bpi^0)
\tensor V(\bpi^0_\zeta)\cdots\tensor V(\bpi^0_{\zeta^{m-2}})
$$
which maps~$v_{\bpi}$ to the corresponding permuted tensor product
of the~$v_{\bpi^0_{\zeta^k}}$. Define~$\eta:V(\bpi)\to V(\bpi)$ by
$\eta:=(\phi_\zeta^*)^{\tensor m}\circ \tau$. Then, for all~$x\in\bu_q$
homogeneous of degree~$k$ and for all~$v\in V(\bpi)$ 
we have~$\eta(xv)=\zeta^{-k} x\eta(v)$ (cf.~\cite[Lemma~2.6]{CG}). 
In particular,
since~$\eta(v_{\bpi})=v_{\bpi}$, we conclude that
$$
V(\bpi)=\bigoplus_{k=0}^{m-1} V(\bpi)^{(k)},\qquad
\text{where~$V(\bpi)^{(k)}:=\{v\in V(\bpi)\,:\, \eta(v)=\zeta^k v\}$}.
$$
Define~$\wh\eta:\wh V(\bpi)\to \wh V(\bpi)$ by~$\wh\eta(v\tensor t^r)=
\zeta^r \eta(v)\tensor t^r$. Then~$\wh\eta\in\End_{\wh\bu_q}\wh V(\bpi)$
(cf.~\cite[Lemma~2.7]{CG}). Moreover, by~\cite[Lemma~2.8]{CG}, 
$\wh V(\bpi)$ is a direct
sum of simple $\wh\bu_q$-submodules~$\wh V(\bpi)^{(r)}$, $r=0,\dots,m-1$
which are in turn the eigenspaces of~$\wh\eta$ corresponding to the
eigenvalues~$\zeta^r$. Observe also that~$\wh V(\bpi)^{(r)}$
is spanned by~$v\tensor t^s$, where~$v\in V(\bpi)^{(k)}$, $k=r-s\pmod m$.
By~\cite[Theorem~5]{CG}, all simple integrable admissible
$\wh\bu_q$-modules of level zero are obtained that way.

\subl{ZCB62}
Following~\cite[4.3]{CG}, set, for all~$v\in V(\bpi)$, $r,s\in\bz$
$$
\Pi_s(v):=\frac1m \sum_{j=0}^{m-1} \zeta^{-js} \eta^j(v),\qquad
\wh\Pi_s(v\tensor t^r):=\Pi_{s-r}(v)\tensor t^r.
$$
By~\cite[Lemma~4.3]{CG}, $\Pi_s$ (respectively, $\wh\Pi_s$) is an orthogonal
projector onto~$V(\bpi)^{(s)}$ (respectively, onto~$\wh V(\bpi)^{(s)}$).
Moreover, if~$x\in\bu_q$ is homogeneous of degree~$k$, then
\begin{equation}\lbl{10}
\Pi_s(x v)=\frac1m \sum_{j=0}^{m-1} \zeta^{-j(s+k)} x \eta^j(v)=
x \Pi_{s+k}(v).
\end{equation}
The map~$\wh\Pi_s$ is obviously a homomorphism of~$\wh\bu_q$-modules.

In the reminder of this section we will prove that~$\wh V(\bpi)^{(r)}$
admits a $\zeta$-crystal basis provided that~$V(\bpi^0)$ admits
a crystal basis.

\subl{ZCB70}
Suppose that~$V(\bpi^0)$ is a ``good''~$\bu_q$-module (we refer the
reader to~\cite[Sect.~8]{Ka2} for the precise definition). In particular,
$V(\bpi^0)$ admits a crystal basis~$(L(\bpi^0),B(\bpi^0))$ and
$B(\bpi^0)^{\tensor m}$ is indecomposable as a crystal for all~$m>0$.
It is proved in~\cite[Proposition~5.15]{Ka2} that the 
module~$V(\bsm\varpi_{i;1})$
corresponding to~$\bpi^0=\bsm\varpi_{i;1}=(\pi_1,\dots,\pi_\ell)$, 
where~$\pi_j(u)=\delta_{i,j}(1-u)$, is good.

Let~$z_1,z_2\in\bc^\times$.
Let~$\tau_{z_1,z_2}$ be the isomorphism
$V(\bpi^0_{z_1})\tensor V(\bpi^0_{z_2})\to V(\bpi^0_{z_2})\tensor 
V(\bpi^0_{z_1})$ normalized so that it preserves the tensor
product of highest weight vectors. By~\cite[Proposition~9.3]{Ka2},
$\tau_{z_1,z_2}$ maps~$\phi_{z_1}^*L(\bpi^0)\tensor_\Ca
\phi_{z_2}^*L(\bpi^0)$ into~$\phi_{z_2}^*L(\bpi^0)\tensor
\phi_{z_1}^*L(\bpi^0)$. Moreover, there is a unique map~$\chi:
B(\bpi^0)^{\tensor 2}\to\bz$ such that
$$
\tau_{z_1,z_2}(b_1\tensor b_2)=(z_1/z_2)^{\chi(b_1\tensor b_2)}
b_1\tensor b_2\pmod{q(\phi_{z_1}^*L(\bpi^0)\tensor_\Ca
\phi_{z_2}^*L(\bpi^0))}.
$$
and~$\chi(b_{\bpi^0}\tensor b_{\bpi^0})=0$ where~$b_{\bpi^0}\in B(\bpi^0)$
is the $l$-highest weight vector.
\begin{lem}
Let~$b_1,b_2\in B(\bpi^0)$ and
suppose that~$\tilde f_i(b_1\tensor b_2)\not=0$. Then
$$
\chi(\tilde f_i(b_1\tensor b_2))=
\begin{cases} \chi(b_1\tensor b_2)+\delta_{i,0},\qquad
&\varphi_i(b_1)>\varepsilon_i(b_2)\\
\chi(b_1\tensor b_2)-\delta_{i,0},& \varphi_i(b_1)\le 
\varepsilon_i(b_2).
\end{cases}
$$
Similarly, if~$\tilde e_i(b_1\tensor b_2)\not=0$, then
$$
\chi(\tilde e_i(b_1\tensor b_2))=
\begin{cases} \chi(b_1\tensor b_2)-\delta_{i,0},\qquad
&\varphi_i(b_1)\ge\varepsilon_i(b_2)\\
\chi(b_1\tensor b_2)+\delta_{i,0},& \varphi_i(b_1)< 
\varepsilon_i(b_2).
\end{cases}
$$
\end{lem}
\begin{pf}
Observe that~$\tilde f_i$ commutes with~$\tau_{z_1,z_2}$. Indeed,
given~$u\in V(\bpi^0_{z_1})\tensor V(\bpi^0_{z_1})$, write, as 
in~\eqref{ZCB10.10}, $u=\sum_{s} F_i^{(s)} u_s$. Since~$\tau_{z_1,z_2}$
is an isomorphism of~$\bu_q$-modules, 
$\tau_{z_1,z_2}(\tilde f_i u)=\sum_s F_i^{(s+1)} \tau_{z_1,z_2}(u_s)$.
On the other hand, $\tau_{z_1,z_2}$ commutes with~$E_i$, $K_i^{\pm1}$, hence
$\tau_{z_1,z_2}(u_s)$ has the same weight as~$u_s$ and is annihilated
by~$E_i$. It follows that 
$\tau_{z_1,z_2}(u)=\sum_s F_i^{(s)} \tau_{z_1,z_2}(u_s)$ 
is the unique decomposition of the form~\eqref{ZCB10.10}. Therefore,
$\tilde f_i \tau_{z_1,z_2}(u)=\sum_s F_i^{(s+1)} \tau_{z_1,z_2}(u_s)=
\tau_{z_1,z_2}(\tilde f_i u)$.

It is sufficient to prove the formula for $\chi(\tilde f_i(b_1
\tensor b_2))$ since the formula for~$\chi(\tilde e_i(b_1\tensor b_2))$
follows from that one by the properties of crystals.
Suppose that~$\varphi_i(b_1)>\varepsilon_i(b_2)$, the other case being
similar. Then
$\tilde f_i(b_1\tensor b_2)=z_1^{-\delta_{i,0}}\tilde f_i b_1\tensor b_2$
by~\lemref{ZCB15} and~\propref{ZCB20}. Therefore,
$$
\tau_{z_1,z_2}(\tilde f_i(b_1\tensor b_2))=z_1^{-\delta_{i,0}}
(z_1/z_2)^{\chi(\tilde f_i b_1\tensor b_2)} \tilde f_i b_1\tensor b_2.
$$
On the other hand, 
$$
\tilde f_i(\tau_{z_1,z_2}(b_1\tensor b_2))=(z_1/z_2)^{\chi(b_1\tensor b_2)}
z_2^{-\delta_{i,0}} \tilde f_i b_1\tensor b_2. 
$$
Since~$\tilde f_i$ commutes with~$\tau_{z_1,z_2}$ it follows that
$\chi(\tilde f_i b_1\tensor b_2)=\chi(b_1\tensor b_2)+\delta_{i,0}$.
\end{pf}
The map~$\chi:B(\bpi^0)^{\tensor 2}\to\bz$ is called the energy function.

\subl{ZCB80}
Retain the notations of~\ref{ZCB60}. 
Using the isomorphism~$\tau_{z_1,z_2}$, we can write~$\tau$ as
$$
\tau=\tau^{(0)}\circ\cdots\circ \tau^{(m-2)},
$$
where
$$
\tau^{(k)}:=\id^{\tensor k}\tensor \tau_{\zeta^k,\zeta^{m-1}}\tensor
\id^{\tensor m-k-2}.
$$
Take some~$b_1,\dots,b_m\in B(\bpi^0)$ and consider~$b=b_1\tensor
\cdots \tensor b_m\in B(\bpi)$. Then
$$
\tau(b)=\zeta^{\Maj_\chi(b)} b,
$$
where~$\chi:B(\bpi^0)^{\tensor 2}\to\bz$ is the energy function
and
$$
\Maj_\chi(b)=\sum_{r=1}^{m-1} r\chi(b_r\tensor b_{r+1})
$$
is the generalised major index of MacMahon. Indeed, in the
case of~$\lie g$ of type~$A_\ell$ and~$\bpi^0=\bsm\varpi_{1;1}$,
there exists a total order on~$B(\bpi^0)$ such that, for all~$b,b'\in
B(\bpi^0)$, $\chi(b\tensor b')=0$
if~$b\ge b'$ whilst~$\chi(b\tensor b')=1$ if~$b<b'$ (cf.~\cite{G}).
Thus, in that case~$\Maj_\chi(b)$ is just the usual major index
of MacMahon for a word in a monoid over a completely ordered alphabet.
\begin{lem}
Let~$V(\bpi)^{(k)}$ be the eigenspace of~$\eta$ corresponding to the
eigenvalue~$\zeta^k$. Set~$L(\bpi)^{(k)}:=L(\bpi)\cap V(\bpi)^{(k)}$,
$B(\bpi)^{(k)}:=\{b\in B(\bpi)\,:\, \Maj_\chi(b)=k\pmod m\}$. Then
\begin{enumerit}
\item $L(\bpi)^{(k)}$ is a free $\Ca$-module, $V(\bpi)^{(k)}=L(\bpi)^{(k)}
\tensor_{\Ca}\bc(q)$ and $B(\bpi)^{(k)}$ is a basis of
the $\bc$-vector space~$L(\bpi)^{(k)}/q L(\bpi)^{(k)}$. 
\item Let~$u\in L(\bpi)^{(k)}$ and write~$u=\sum_{s} F_i^{(s)}u_s$ as
in~\eqref{ZCB10.10}. Then~$u_s\in L(\bpi)^{(k-s\delta_{i,0})}$,
$\tilde e_i u\in L(\bpi)^{(k-\delta_{i,0})}$ and $\tilde f_i u\in
L(\bpi)^{(k+\delta_{i,0})}$.
\item Suppose that~$b\in B(\bpi)^{(k)}$. Then 
$$
\tilde e_i b\in \zeta^{\bz
\delta_{i,0}} B(\bpi)^{(k-\delta_{i,0})}\cup\{0\},\quad\tilde f_i b\in
\zeta^{\bz\delta_{i,0}} B(\bpi)^{(k+\delta_{i,0})}\cup\{0\}.
$$
\end{enumerit}
\end{lem}
\begin{pf*}
Take~$u\in L(\bpi)$ such that~$u=b\pmod{q L(\bpi)}$. Since~$L(\bpi)$ is
a free module and~$B(\bpi)$ is a basis of~$L(\bpi)/qL(\bpi)$, such~$u$
generate $L(\bpi)$ as an~$\Ca$-module by Nakayama's Lemma. Then, 
since~$\eta$ maps $L(\bpi)$ into itself,
$$
\Pi_s(u)=\frac1m \sum_{r=0}^{m-1} \zeta^{r(\Maj_\chi(b)-s)} b\pmod{q L(\bpi)}.
$$
It follows that~$\Pi_s(u)=b\pmod{qL(\bpi)}$ if~$s=\Maj_\chi(b)\pmod m$ 
whilst~$\Pi_s(u)=0\pmod{qL(\bpi)}$ otherwise. 

Since~$\Pi_k$ is an orthogonal projector onto~$V(\bpi)^{(k)}$ and
maps~$L(\bpi)$ into itself, it follows that~$L(\bpi)^{(k)}=\Pi_k(L(\bpi))$.
Then~$B(\bpi)^{(k)}$ is a basis of~$L(\bpi)^{(k)}/q L(\bpi)^{(k)}$.
Indeed, elements of~$B(\bpi)^{(k)}$ are contained 
in~$L(\bpi)^{(k)}/q L(\bpi)^{(k)}$ by the above and are linearly independent,
whence~$\dim_\bc L(\bpi)^{(k)}/q L(\bpi)^{(k)}\ge \# B(\bpi)^{(k)}$. Yet,
$\# B(\bpi)=\sum_{k=0}^{m-1}\# B(\bpi)^{(k)}\le 
\sum_{k=0}^{m-1} \dim_\bc L(\bpi)^{(k)}/q L(\bpi)^{(k)}=\dim_\bc
L(\bpi)/q L(\bpi)=\# B(\bpi)$. It follows that
$\dim_\bc L(\bpi)^{(k)}/q L(\bpi)^{(k)}=\# B(\bpi)^{(k)}$. Then
$L(\bpi)^{(k)}$ is generated by the~$\Pi_k(u)$, $u=b\pmod{q L(\bpi)}$,
with~$\Maj_\chi(b)=k$ by Nakayama's Lemma.

For the second part, suppose that~$u\in L(\bpi)^{(k)}$. 
Then by~\eqref{ZCB62.10},
$$
u=\Pi_k(u)=\sum_{s} \Pi_k( F_i^{(s)} u_s)=\sum_{s} F_i^{(s)} 
\Pi_{k-s\delta_{i,0}}(u_s).
$$
Since $K_i$ commutes with the~$\Pi_r$ and
$E_i\Pi_{r}(u_s)=\Pi_{r-\delta_{i,0}}(E_iu_r)=0$ it follows that
$\Pi_{k-s\delta_{i,0}}(u_s)$ is of the same weight as~$u_s$
and is annihilated by~$E_i$. Then $u_s=\Pi_{k-s\delta_{i,0}}(u_s)$
by the uniqueness of the decomposition~\eqref{ZCB10.10}. Furthermore,
$$
\Pi_r(\tilde e_i u)=\sum_{s} \Pi_r (F_i^{(s-1)}u_s)=
\sum_{s} F_i^{(s-1)} \Pi_{r-(s-1)\delta_{i,0}}(u_s).
$$
It remains to observe that~$\Pi_{r-(s-1)\delta_{i,0}}(u_s)=0$ unless~$r=k-
\delta_{i,0}$. The proof for~$\tilde f_i$ is similar.

The last part follows immediately from~(i), (ii) and the properties of the
$z$-crystal basis. However, we prefer to present a direct proof since 
it involves a property of~$\Maj_\chi$ which we will need later. 
Evidently, it is enough to prove the statement for~$\tilde f_i$.
Write~$b=b_1\tensor \cdots \tensor b_m$, $b_i\in B(\bpi^0)$ and
suppose that~$\tilde f_i b\not=0$. 
Then~$\tilde f_i b=\zeta^{-s+1} b_1\tensor \cdots\tensor
\tilde f_i b_s\tensor\cdots \tensor b_{m}$ for some~$1\le s\le m$. Suppose
first that~$1\le s<m$. Then
\begin{align*}
\Maj_\chi(\tilde f_i b)-\Maj_\chi(b)&=
(s-1)(\chi(b_{s-1}\tensor \tilde f_i b_s)-
\chi(b_{s-1}\tensor b_s))\\&+s(\chi(\tilde f_i b_s\tensor b_{s+1})
-\chi(b_s\tensor b_{s+1}))\\
&=-(s-1)\delta_{i,0}+s\delta_{i,0}=\delta_{i,0},
\intertext{
where we used~\lemref{ZCB70}. Finally, if~$s=m$, then}
\Maj_\chi(\tilde f_i b)-\Maj_\chi(b)&=(m-1)(\chi(b_{m-1}
\tensor\tilde f_i b_m)-\chi(b_{m-1}\tensor b_m))\\&=-(m-1)\delta_{i,0}=
\delta_{i,0}\pmod m.\tag*\qed
\end{align*}
\end{pf*}

\subl{ZCB90}
Retain the notations of~\ref{ZCB30}, \ref{ZCB60} and~\ref{ZCB80}.
\begin{thm*}
Suppose that~$V(\bpi^0)$ is a good module and let~$(L(\bpi^0),
B(\bpi^0))$ be its crystal basis. Set~$\bpi=
\bpi^0 \bpi^0_\zeta\cdots \bpi^0_{\zeta^{m-1}}$, where $\zeta$
is an $m$th primitive root of unity, and define~$L(\bpi)$,
$B(\bpi)$ as in~\ref{ZCB30}. 
The simple submodule~$\wh V(\bpi)^{(k)}$, $k=0,\dots,m-1$ of~$\wh
V(\bpi)$ admits a $\zeta$-crystal base~$(\wh L(\bpi)^{(k)}, 
\wh B(\bpi)^{(k)})$, where
\begin{align*}
&\wh L(\bpi)^{(k)}=\wh\Pi_k \wh L(\bpi)=\bigoplus_{\substack{r\in\bz,\,
0\le s\le m-1\\r+s=k\pmod m}}\mskip-15mu
L(\bpi)^{(s)}\tensor t^{r},\\
&\wh B(\bpi)^{(k)}=\{ b\tensor t^r\,:\, b\in B(\bpi)^{(s)},\,r\in\bz,\,
r+s=k\pmod m\}
\end{align*}
\end{thm*}
\begin{pf}
This follows immediately from~\lemref{ZCB40} and~\lemref{ZCB80}.
\end{pf}
Our~\thmref{thm1} is a particular case of the above statement
since, as shown in~\cite{Ka2}, the module~$V(\bpi^0)$
with~$\bpi^0=\bsm\varpi_{i;1}$ satisfies all the required conditions
and the corresponding~$\bpi$ obviously coincides with~$\bsm\varpi_{i;m}$.

\section{Path model for $z$-crystal bases of quantum loop modules}

In the present section we will construct a combinatorial model,
in the framework of Littelmann's path crystal, of $z$-crystal bases
of simple components of quantum loop modules of fundamental type.
The necessary facts about Littelmann's path crystal will be reviewed
as the need arises. Throughout this section we identify~$P$ with~$\wh P/
\bz\delta$.

\subl{PM10}
Given~$a,b\in\bq$, $a<b$, set~$[a,b]:=\{x\in\bq\,|\, a\le x\le b\}$.
Let~$\mathbb P$ (respectively, $\wh{\mathbb P}$) be the set
of piece-wise linear continuous paths in~$P\tensor_\bz\bq$ (respectively,
in~$\wh P\tensor_\bz\bq$) starting at zero and terminating at
an element of~$P$ (respectively, $\wh P$). In other words,
$\pi\in\mathbb P$ (respectively $\hatbp$) is a piece-wise linear continuous
map of~$[0,1]$ into~$P\tensor_\bz\bq$ 
(respectively, into $\wh P\tensor_\bz\bq$)
such that~$\pi(0)=0$ and~$\pi(1)\in P$ (respectively, $\wh P$).
We consider two paths as identical if they coincide up to a continuous
piece-wise linear non-decreasing reparametrisation.

After Littelmann (cf.~\cite{Li94,Li95}) one can introduce a structure
of a normal crystal on~$\mathbb P$ or on~$\wh{\mathbb P}$ in the following
way. Given~$\pi\in\mathbb P$ or~$\hatbp$ and~$i\in\wh I$, set
$h^i_\pi(\tau)=-\<\alpha_i^\vee,\pi(\tau)>$, $\tau\in[0,1]$. 
Let~$\varepsilon_i(\pi)$ be the maximal integral value attained
by~$h^i_\pi$ on~$[0,1]$. Furthermore, set~$e_+^i(\pi)=\min\{\tau\in[0,1]\,:\,
h^i_\pi(\tau)=\varepsilon_i(\pi)\}$. If~$\varepsilon_i(\pi)=0$ then set~$e_i\pi=0$.
Otherwise, set~$e_-^i(\pi)=\max\{\tau\in[0,e_i^+(\tau)]\,:\,
h^i_\pi(\tau)=\varepsilon_i(\pi)-1\}$ and define
$$
(e_i\pi)(\tau)=\begin{cases}
\pi(\tau),\qquad&\tau\in[0,e_-^i(\pi)]\\
\pi(e_-^i(\pi))+s_i(\pi(\tau)-\pi(e_-^i(\pi))),&\tau\in[e_-^i(\pi),
e_+^i(\pi)]\\
\pi(\tau)+\alpha_i,&\tau\in[e_+^i(\pi),1],
\end{cases}
$$
where~$s_i$ acts point-wise. Similarly, in order to define~$f_i$,
let~$f_+^i(\pi)=\max\{\tau\in[0,1]\,:\, h^i_\pi(\tau)=\varepsilon_i(\tau)\}$.
If~$f_+^i(\pi)=1$, set~$f_i\pi=0$. Otherwise, set~$f_-^i(\pi)=
\min\{\tau\in[f_+^i(\pi),1]\,:\,h^i_\pi(\tau)=\varepsilon_i(\pi)-1\}$
and define
$$
(f_i\pi)(\tau)=\begin{cases}
\pi(\tau),\qquad&\tau\in[0,f_+^i(\pi)]\\
\pi(f_+^i(\pi))+s_i(\pi(\tau)-\pi(f_+^i(\pi))),&\tau\in[f_+^i(\pi),
f_-^i(\pi)]\\
\pi(\tau)-\alpha_i,&\tau\in[f_-^i(\pi),1],
\end{cases}
$$
Finally, $\wt\pi$ is defined as the endpoint~$\pi(1)$ of~$\pi$.
\begin{rem}
As in~\cite{G}, we use the definition of crystal operations on~$\mathbb P$ 
given in~\cite[6.4.4]{JB} which differs by the sign of~$h^i_\pi$ from
the definition in~\cite[1.2]{Li94}.
That choice is more convenient for us since it makes the comparison with
Kashiwara's tensor product easier.
\end{rem}

\subl{PM40}
Following~\cite[Theorem~8.1]{Li95}, one can introduce
an action of the Weyl group~$\wh W$ on~$\mbp$ and~$\hatbp$.
Namely, given~$\pi\in\mathbb P$ or~$\hatbp$, set
$$
s_i \pi=\begin{cases}
f_i^{\<\alpha_i^\vee,\pi(1)>}\pi,\qquad &\<\alpha_i^\vee,\pi(1)> \ge 0,\\
e_i^{-\<\alpha_i^\vee,\pi(1)>}\pi,&\<\alpha_i^\vee,\pi(1)> \le 0.
\end{cases}
$$
Given~$\lambda\in P$ or~$\wh P$, denote by~$\pi_\lambda$
the linear path~$\tau\mapsto \tau\lambda$. One can easily see
from the definitions in~\ref{PM10} that~$\varepsilon_i(\pi_\lambda)=\max\{0,
-\<\alpha_i^\vee,\lambda>\}$ and $\varphi_i(\pi_\lambda)=
\max\{0,\<\alpha_i^\vee,\lambda>\}$. 
\begin{lem}
For all~$\lambda\in P$ or~$\wh P$, $s_i\pi_\lambda=\pi_{s_i\lambda}$.
In particular, if $B$ is a subcrystal of~$\mathbb P$ or~$\hatbp$ 
and $\pi_\lambda\in B$ for some~$\lambda\in P$ or~$\wh P$ then
$\pi_{w\lambda}\in B$ for all~$w\in\wh W$.
\end{lem}
\begin{pf}
The second assertion is an immediate corollary of the first one which
in turn follows from the formulae
\begin{alignat*}{2}
&f_i^n \pi_\lambda =\begin{cases}
s_i\lambda\tau,\qquad&\tau\in\big[0,\frac{n}{\<\alpha_i^\vee,\lambda>}\big]
\\
\lambda\tau-n\alpha_i,&\tau\in\big[\frac{n}{\<\alpha_i^\vee,\lambda>},1\big]
\end{cases},&\qquad&0< n\le \<\alpha_i^\vee,\lambda>\\
&e_i^n \pi_\lambda =\begin{cases}
\lambda\tau,\qquad&\tau\in\big[0,1-\frac{n}{|\<\alpha_i^\vee,\lambda>|}
\big]\\
s_i\lambda\tau+(|\<\alpha_i^\vee,\lambda>|-n)\alpha_i,
&\tau\in\big[1-\frac n{|\<\alpha_i^\vee,\lambda>|},1\big]
\end{cases},&& 0< n\le -\<\alpha_i^\vee,\lambda>.
\end{alignat*}
These can be deduced easily from the formulae in~\ref{PM10} 
by induction on~$n$.
\end{pf}

\subl{PM50}
Given~$\lambda\in P$ or~$\wh P$ and~$\mu,\nu\in \wh W\lambda$, write,
following~\cite{Li95}, 
$\nu\ge\mu$ if there exist a sequence~$\{\nu_0=\nu,\nu_1,\dots,\nu_s=\mu\}$,
$\nu_i\in P$ or~$\wh P$ 
and positive real roots~$\beta_1,\dots,\beta_s$ of~$\hatg$
such that
$$
\nu_i=s_{\beta_i}(\nu_{i-1}),\qquad\<\beta_i^\vee,\nu_{i-1}><0,\qquad
i=1,\dots,s.
$$
If~$\nu\ge\mu$, let~$\dist(\nu,\mu)$ be the maximal length of such a 
sequence.

Let~$\bsm\nu=\{\nu_1,\dots,\nu_r\}$ be a sequence of
elements of~$\wh W\lambda$ and~$\bsm a=\{a_0=0<a_1<\cdots<a_r=1\}$ be
a sequence of rational numbers. Denote by~$\pi_{\bsm\nu,\bsm a}$
the piece-wise linear path
\begin{equation}\lbl{10}
\pi_{\bsm\nu,\bsm a}(\tau)=\sum_{i=1}^{j-1}(a_i-a_{i-1})\nu_i+
(t-a_{j-1})\nu_j,\qquad \tau\in[a_{j-1},a_j].
\end{equation}
In other words, it is a concatenation of straight lines joining~$\lambda_{j-1}$
and~$\lambda_j$, $j=0,\dots,r$, where~$\lambda_j=\sum_{i=1}^j
(a_i-a_{i-1})\nu_i$.
\begin{defn}[{\cite{Li95}}]
Fix~$\lambda\in P$ or~$\wh P$. 
A path of the form~$\pi_{\bsm\nu,\bsm a}$, where~$\bsm\nu=\{\nu_1
\ge\cdots\ge\nu_r\}$, $\nu_i\in \wh W\lambda$ and~$\bsm a=\{a_0=0<a_1<
\cdots<a_r=1\}$ is called a {\em Lakshmibai-Seshadri
(LS) path of class~$\lambda$} 
if, for all~$1\le i\le r-1$, either~$\nu_i=\nu_{i+1}$
or there exists a sequence~$\lambda_{0,i}=\nu_i>\lambda_{1,i}>\cdots>
\lambda_{s,i}=\nu_{i+1}$, $\lambda_{j,i}\in\wh W\lambda$ such that
$$
\lambda_{j,i} = s_{\beta_{j,i}}(\lambda_{j-1,i}),\quad
a_i\<\beta_{j,i}^\vee,\lambda_{j-1,i}>\in-\bn,\quad
\dist(\lambda_{j-1,i},\lambda_{j,i})=1,
$$
for some real positive roots~$\beta_{j,i}$.
\end{defn}
It is known (cf.~\cite[Lemma~4.5]{Li95}) that an LS-path~$\pi=\pi_{\bsm\nu,
\bsm a}$ 
of class~$\lambda$
is an element of~$\mathbb P$ or~$\hatbp$ and has the integrality property, that
is, the maximal value attained by the function~$h^i_\pi$ on~$[0,1]$
is an integer for all~$i\in\wh I$. Moreover, by~\cite[Lemma~4.5]{Li95}
all local maxima of~$h^i_\pi$ are integers.

\subl{PM55}
Given a collection~$\pi_1,\dots,\pi_k$ of paths in~$\mathbb P$ or~$\hatbp$, 
define their concatenation
$$
(\pi_1\tensor\cdots\tensor \pi_k)(\tau)=
\sum_{1\le s<j} \pi_s(1)+\pi((\tau-\sigma_{j-1})/(\sigma_j-\sigma_{j-1})),
\qquad \tau\in[\sigma_{j-1},\sigma_j]
$$
for some~$0=\sigma_0<\sigma_1<\cdots<\sigma_{k-1}<\sigma_k=1$,
$\sigma_j\in\bq$. 
This definition does not depend on the~$\sigma_j$, up to a reparametrisation.
By~\cite[2.6]{Li95}, the concatenation of
paths satisfies Kashiwara's tensor product rule.

Let~$\pi\in\mathbb P$ or~$\hatbp$ be an LS-path and define the path~$n\pi$
by~$(n\pi)(\tau)=n \pi(\tau)$, $\tau\in[0,1]$. Evidently,
$\varepsilon_i(n\pi)=n\varepsilon_i(\pi)$, $\varphi_i(n\pi)=n\varphi_i(\pi)$,
$i\in\wh I$ and $\wt n\pi=n\wt \pi$. Let~$S_n$ be the map~$\pi\mapsto n\pi$.
Then by~\cite[Lemma~2.4]{Li95}, $S_n(e_i\pi)=e_i^n S_n(\pi)$ and~$
S_n(f_i\pi)=f_i^n S_n(\pi)$, $i\in\wh I$. Observe that, for a linear
path~$\pi_\lambda$, $S_n(\pi_\lambda)=\pi_{n\lambda}=\pi_\lambda^{\tensor n}$.

\subl{PM60}
Fix~$\lambda\in P$ (respectively, $\lambda\in\wh P$) 
and let~$B(\lambda,\mbp)$ (respectively, $B(\lambda,\hatbp)$) 
be the subcrystal of~$\mathbb P$ (respectively, $\hatbp$)
generated over the monoid $\Cm$ (cf.~\ref{P50}) by the linear
path~$\pi_\lambda$. Henceforth we write~$B(\lambda)$ for~$B(\lambda,\mbp)$.
Then by~\cite[Corollary~2 of Proposition~4.7]{Li95}
all elements of~$B(\lambda)$ or~$B(\lambda,\hatbp)$ are LS-paths of
class~$\lambda$. 
Suppose further that~$B(\lambda)$ is a finite set.
Then there exists~$N\in\bn^+$ and~$\bsm a=\{0=a_0<a_1<\cdots<a_N=1\}$
such that every element of~$B(\lambda)$
can be represented as~$\pi_{\bsm\nu,\bsm a}$ for some
sequence of weights~$\bsm\nu=
\{\nu_1\ge \cdots\ge \nu_N\}$. Observe that~$\pi_{\bsm\nu',\bsm a'}(\tau)=
\pi_{\bsm\nu,\bsm a}(\tau)$, for all~$\tau\in[0,1]$, 
where
\begin{alignat*}{3}
&a'_j = a_j,&\quad & \nu'_j = \nu_j,&\quad & j=0,\dots,r\\
&a'_{r+1}=x,&& \nu'_{r+1}=\nu_{r+1},\\
&a'_j = a_{j-1},&& \nu'_j = \nu_{j-1},&& j=r+2,\dots,N+1
\end{alignat*}
for any~$0\le r<N$ and for any rational~$x$, $a_r < x < a_{r+1}$.
Therefore we may assume, without loss of generality,
that~$a_j=j/N$ and in that case we omit~$\bsm a$.
\begin{lem}
\begin{enumerit}
\item $S_N(\pi_{\bsm\nu})=\pi_{\nu_1}\tensor\cdots\tensor \pi_{\nu_N}$.
In particular, $S_N$ can be viewed as an injective map~$B(\lambda)
\to B(\lambda)^{\tensor N}$.
\item If~$e_i\pi_{\bsm\nu}\not=0$ then 
$e_i\pi_{\bsm\nu}=\pi_{\bsm\nu'}$ with
$$
\bsm\nu'=\{\nu_1,\dots,\nu_k,s_i(\nu_{k+1}),\dots,
s_i(\nu_l),\nu_{l+1},\dots,\nu_N\}
$$ where~$k=Ne^i_-(\pi_{\bsm\nu})$
and~$l=Ne^i_+(\pi_{\bsm\nu})$. Similarly, 
if~$f_i\pi_{\bsm\nu}\not=0$ then $f_i\pi_{\bsm\nu}
=\pi_{\bsm\nu''}$ with
$$
\bsm\nu''=\{\nu_1,\dots,\nu_r,s_i(\nu_{r+1}),\dots,
s_i(\nu_s),\nu_{s+1},\dots,\nu_N\},
$$
where~$r=Nf^i_+(\pi_{\bsm\nu})$ and~$s=Nf^i_-(\pi_{\bsm\nu})$.
\item Let~$k$, $l$, $r$ and~$s$ be as above. Then
$$
\sum_{j=k+1}^l \<\alpha_i^\vee,\nu_j>=-N,\qquad
\sum_{j=r+1}^s \<\alpha_i^\vee,\nu_j>=N.
$$
\end{enumerit}
\end{lem}
\begin{pf}
Let~$\pi=\pi_{\bsm\nu}$. 
The first part follows immediately from~\eqref{PM50.10} and~\ref{PM55}
with~$\sigma_j=j/N$, $j=0,\dots,N$. In order to prove~(iii) observe
that, by the choice of~$\bsm a$, there exist~$0\le k<l\le N$ such that
$k/N=e^i_-(\pi)$ whilst~$l/N=e^i_+(\pi)$. Then 
$$
1=h^i_\pi(e^i_+(\pi))-h^i_\pi(e^i_-(\pi))=
-\frac1N\sum_{r=k+1}^l \<\alpha_i^\vee,\nu_r>,
$$
by~\eqref{PM50.10}. The second formula in~(iii) is proved in the same
way. 
Furthermore, for~$\tau\in[(j-1)/N,j/N]$ one has by~\eqref{PM50.10}
\begin{alignat*}{2}
(e_i \pi)(\tau)&=\frac1N\sum_{r=1}^{j-1} \nu_r+
(\tau-(j-1)/N)\nu_j,&\,& 1\le j\le k\\
(e_i \pi)(\tau)&=\frac1N\sum_{r=1}^{j-1} s_i \nu_r
+(\tau-(j-1)/N)s_i\nu_j+\pi(k/N)-s_i\pi(k/N)\\
&=\frac1N\sum_{r=1}^k \nu_r+\frac1N\sum_{r=k+1}^j
s_i\nu_r +(\tau-(j-1)/N)s_i\nu_j,&&
k\le j\le l\\
(e_i \pi)(\tau)&=\frac1N\sum_{r=1}^{j-1} \nu_r
+(\tau-(j-1)/N)\nu_j+\alpha_i,&&l\le j\le N.
\end{alignat*}
It is now obvious that $(e_i\pi)(\tau)=\pi_{\bsm\nu'}(\tau)$, $0\le \tau
\le l/N$. Finally, observe that~$
\alpha_i=-\frac1N\sum_{r=k+1}^l
\<\alpha_i^\vee,\nu_r>\alpha_i=
\frac1N\sum_{r=k+1}^l (s_i\nu_r-\nu_r)$ by~(iii). Thus, we can
write for~$\tau\in[(j-1)/N,j/N]$ and~$l\le j \le N$
$$
(e_i\pi)(\tau)=\frac1N\Big(\sum_{\substack{
1\le r\le k,\\
l< r\le N}}\nu_r+\sum_{k< r\le l} s_i\nu_r\Big)+
(\tau-(j-1)/N)\nu_j=\pi_{\bsm\nu'}(\tau).
$$
The second formula in~(ii) for the action of~$f_i$ is proved in a similar way.
\end{pf}

\subl{PM65}
Let~$\xi:\wh P\to P$ be the canonical projection. Define the map
$\Xi:\hatbp\to\mbp$ by~$(\Xi \pi)(\tau)=\xi(\pi(\tau))$, for all~$\tau\in
[0,1]$.
\begin{lem}
The map~$\Xi$ is a morphism of crystals and~$\Xi(B(\lambda,\hatbp))=
B(\xi(\lambda))$.
\end{lem}
\begin{pf}
Since~$\<\alpha_i^\vee,\delta>=0$ for all~$i\in\wh I$, we conclude
that~$h^i_{\Xi\pi}(\tau)=h^i_{\pi}(\tau)$ for all~$i\in\wh I$ and
for all~$\tau\in[0,1]$. For the same reason, one has
$\xi(s_i(\pi(\tau)))=s_i((\Xi \pi)(\tau))$. It is now obvious 
that~$\Xi$ commutes with the operators defined in~\ref{PM10}.
In order to prove the second assertion, observe that~$\Xi\pi_\lambda=
\pi_{\xi(\lambda)}$.
\end{pf}

\subl{PM70}
It was shown in~\cite[Theorem~1.1]{NS} that 
$B(\varpi_i,\hatbp)$ is isomorphic to the crystal basis~$\wh B(
\bsm\varpi_{i;1})$ of the
simple integrable $\wh\bu_q$-module $\wh V(\bsm\varpi_{i;1})$.
Thus, for all~$\pi\in B(\varpi_i,\hatbp)$ there exists~$x\in\Cm$
such that~$\pi=x\pi_{\varpi_i}$. Following~\cite[5.3]{NS},
one introduces a translation operator~$z$ on~$B(\varpi_i,\hatbp)$
by~$z\pi=x\pi_{\varpi_i+\delta}$. It follows that~$(z\pi)(\tau)=
\pi(\tau)+\tau\delta$ for all~$\tau\in[0,1]$.
By~\cite[Proposition~5.8]{NS},
$B(\varpi_i,\hatbp)/\!\sim$, where~$\sim$ is an equivalence relation
on~$B(\varpi_i,\hatbp)$ defined by~$\pi\sim \pi'$ if and only if~$\pi=
z^k \pi'$ for some~$k\in\bz$, is a crystal isomorphic to the
crystal basis~$B(\bsm\varpi_{i;1})$ of the finite dimensional simple
$\bu_q$-module $V(\bsm\varpi_{i;1})$. Evidently, 
$B(\varpi_i,\hatbp)/\!\!\sim$ is isomorphic to~$\Xi(B(\varpi_i,\hatbp))=
B(\varpi_i)\subset \mbp$ in the notations of~\ref{PM65}.

Thus, $B(\varpi_i)$ is isomorphic to~$B(\bsm\varpi_{i;1})$ as a crystal.
In particular, $B(\varpi_i)$ is finite, $B(\varpi_i)^{\tensor m}$
is indecomposable as a crystal for all~$m>0$ (hence~$B(\varpi_i)^{\tensor m}
=B(m\varpi_i)$) and there exists a 
unique map~$\chi:B(\varpi_i)^{\tensor 2}
\to \bz$ satisfying~$\chi(\pi_{\varpi_i}\tensor \pi_{\varpi_i})=0$
and the properties listed in~\lemref{ZCB70}.

We will now construct an injective map~$\psi$ 
from~$\wh{B(\varpi_i)^{\tensor m}}$ into~$\hatbp$.
Given an arbitrary collection of weights~$\bsm\lambda=
\{\lambda_0,\dots,\lambda_K\}$, where~$\lambda_0=0$ and~$\lambda_j
\in P\tensor_\bz\bq$ or~$\wh P\tensor_\bz\bq$, 
and a collection of rational numbers~$\bsm a=\{a_0=0<a_1<\cdots<a_K=1\}$
denote by~$p_{\bsm\lambda,\bsm a}$ the path
$$
p_{\bsm\lambda,\bsm a}(\tau)=\lambda_{j-1}+
(\lambda_j-\lambda_{j-1})\left(\frac{\tau-a_{j-1}}{a_j-a_{j-1}}\right),
\qquad \tau\in[a_{j-1},a_j].
$$
In other words, $p_{\bsm\lambda,\bsm a}$ is a concatenation of straight
lines joining~$\lambda_{j-1}$ with~$\lambda_j$, $j=1,\dots,K$.
As before, we omit~$\bsm a$ if~$a_j=j/K$.

Let~$b=b_1\tensor\cdots\tensor b_m$ be an element of~$B(\varpi_i)^{\tensor m}$
and choose $N\in\bn^+$ as in~\ref{PM60}. Actually, it is sufficient
to take as~$N$ the least common multiple of coefficients of~$\alpha_i^\vee$
in all co-roots of~$\lie g$.
Then~$b_k=\pi_{\bsm\nu^{(k)}}$
with~$\bsm\nu^{(k)}=\{ w^{(k)}_1\varpi_i\ge \cdots\ge
w^{(k)}_N\varpi_i\}$ for some~$w^{(k)}_j\in \wh W$.
Set~$\nu_{rN+s}= w^{(r+1)}_s \varpi_i$, $r=0,\dots,m-1$, $s=1,\dots,N$.
The element~$b$ then corresponds to the linear path~$p_{\bsm\lambda}
\in \mathbb P$ 
with~$\bsm\lambda=\{\lambda_0=0,\lambda_1\dots,\lambda_{Nm}\}$ 
where~$\lambda_j=\frac1N\sum_{k=1}^j \nu_k$.
On the other hand, using the map~$S_N$ and the associativity of the
tensor product of crystals, we associate with~$b$ the element
$T_N(b):= \pi_{1}\tensor\cdots\tensor \pi_{{Nm}}\in B(\varpi_i)^{\tensor Nm}$,
where~$\pi_r:=\pi_{\nu_r}$.
Set $\Maj_\chi(T_N(b))=
\sum_{r=1}^{Nm-1} r\chi(\pi_{r}\tensor\pi_{{r+1}})$.
This expression is defined since~$\nu_r\in\wh W\varpi_i$ for 
all~$r=1,\dots,Nm$ and so~$\pi_{r}\in B(\varpi_i)$ by~\lemref{PM40}

For all~$n\in\bz$ we associate 
with~$b\tensor t^n\in\wh{B(\varpi_i)^{\tensor m}}$ 
a path~$p_{\wh{\bsm\lambda}(n)}\in\hatbp$, where~$\wh{\bsm\lambda}(n)=
\{\wh\lambda_0=0,\wh\lambda_1,\dots,\wh\lambda_{Nm}\}$,
$\wh\lambda_j=\lambda_j+\kappa_j(b\tensor t^n) 
\delta$ and
\begin{align*}
\kappa_j(b\tensor t^n)
=&\frac{j}{Nm}\Big(\frac1N\Maj_\chi(T_N(b))+n\Big)\\&-
\frac1N \sum_{s=1}^{j-1} s \chi(\pi_{s}\tensor \pi_{{s+1}})
-\frac jN \sum_{s=j}^{Nm-1} \chi(\pi_{s}\tensor \pi_{{s+1}}).
\end{align*}
\begin{prop}
The map~$\psi:\wh{B(\varpi_i)^{\tensor m}}\to\hatbp$ sending
$p_{\bsm\lambda}\tensor t^n$ to~$p_{\wh{\bsm\lambda}(n)}$ is an 
injective morphism of crystals.
\end{prop}
\begin{pf}
The injectivity is obvious.
Let~$b=p_{\bsm\lambda}$.
By the definition of affinisation, $\wt b\tensor t^n=
\wt b+n\delta=\lambda_{Nm}+n\delta$. On the other hand,
$\wt \psi(b\tensor t^n)=\lambda_{Nm}+\kappa_{Nm}(b\tensor t^n)\delta$
and
$$
\kappa_{Nm}(b\tensor t^n)=\frac1N\,\Maj_\chi(T_N(b))+n-
\frac1N\sum_{s=1}^{Nm-1}s\chi(\pi_{s}\tensor \pi_{{s+1}})=
n,
$$
hence~$\wt b\tensor t^n=\wt \psi(b\tensor t^n)$. Furthermore,
observe that by construction~$\Xi p_{\wh{\bsm\lambda}(n)}=p_{\bsm\lambda}$.
Then it follows from~\lemref{PM65} that
$\varepsilon_j$ and hence~$\varphi_j$ commute with~$\psi$ 
for all~$j\in\wh I$. 

Since both~$\wh{B(\varpi_i)^{\tensor m}}$ and~$\hatbp$ are normal,
it remains to prove that~$\psi(e_j(b\tensor t^n))=
\psi(e_j b\tensor t^{n+\delta_{j,0}})=e_j \psi(b\tensor t^n)$ provided
that~$e_j b\not=0$.

By the choice of~$N$, there exist~$0\le k<l\le Nm$ such 
that~$e^j_-(b)=k/Nm$ and~$e^j_+(b)=l/Nm$. As above, we conclude
immediately that~$e^j_-(\psi(b\tensor t^n))=k/Nm$ and~$e^j_+(
\psi(b\tensor t^n))=l/Nm$. 
Set~$p_{\bsm\mu}=e_j b$ and~$p_{\wh{\bsm\lambda}{}'}=e_j\psi(b\tensor
t^n)=e_jp_{\wh{\bsm\lambda}(n)}$. 
We aim to prove that~$\wh{\bsm\mu}(n+
\delta_{j,0})=
\wh{\bsm\lambda}{}'$. 
Using the formulae from~\ref{PM10} and the definition of~$\psi$ we obtain
\begin{align*}
&\wh\mu_s=\begin{cases}
\lambda_s+K_s'\delta,\qquad& 0\le s\le k\\
\lambda_s+\<\alpha_j^\vee,\lambda_k-\lambda_s>\xi(\alpha_j)
+K_s'\delta,&
k\le s\le l\\
\lambda_s+\xi(\alpha_j)+K_s'\delta,& l\le s\le Nm,
\end{cases}
\\
\intertext{where~$K_s':=\kappa_s(e_jb\tensor t^{n+\delta_{j,0}})$.
On the other hand,}
&\wh\lambda{}_s'=\begin{cases}
\lambda_s+K_s\delta,\qquad& 0\le s\le k\\
\lambda_s+\<\alpha_j^\vee,\lambda_k-\lambda_s>\alpha_j+K_s\delta,
\hphantom{\xi()}\mskip-2.5mu&
k\le s\le l\\
\lambda_s+\alpha_j+K_s\delta,& l\le s\le Nm,
\end{cases}
\end{align*}
where~$K_s:=\kappa_s(b\tensor t^n)$. Since~$\alpha_j=
\delta_{j,0}\delta+\xi(\alpha_j)$, $j\in\wh I$, the above formulae
imply that~$\wh{\bsm\mu}(n+\delta_{j,0})=\wh{\bsm\lambda}{}'$ if and only if
\begin{alignat}{2}
&K_s'=K_s,&\qquad& 0\le s\le k \lbl{10a}\\
&K_s'=K_s+\<\alpha_j^\vee,\lambda_k-\lambda_s>\delta_{j,0},&& k\le s\le l \lbl{10b}\\
&K_s'=K_s+\delta_{j,0},&& l\le s\le Nm.\lbl{10c}
\end{alignat}

Write~$T_N(e_jb)=\pi_1'\tensor\cdots\tensor \pi_{Nm}'$. By
\lemref{PM60}, $\pi_r'=\pi_{r}$,
$1\le r\le k$ or~$l<r\le Nm$ whilst~$\pi_r'=s_j \pi_{r}=
\pi_{s_j\nu_r}$, $r=k+1,\dots,l$.
Then
\begin{equation}\lbl{20}
\begin{split}
K_s'-K_s=\frac{s}{Nm}\Big(&\frac1N(\Maj_\chi(T_N(e_jb))-
\Maj_\chi(T_N(b)))+\delta_{j,0}\Big)\\
-&\frac1N \sum_{r=1}^{s-1} 
r(\chi(\pi'_r\tensor \pi'_{r+1})-\chi(\pi_{r}\tensor \pi_{{r+1}}))\\
-&\frac sN\sum_{r=s}^{Nm-1}
(\chi(\pi'_r\tensor \pi'_{r+1})-\chi(\pi_{r}\tensor \pi_{{r+1}}))
\end{split}
\end{equation}
Evidently, $\chi(\pi'_r\tensor \pi'_{r+1})=\chi(\pi_r\tensor 
\pi_{r+1})$ for all~$1\le r<k$ and for all~$l<r<Nm$. The crucial
point in our argument is the following
\begin{clm}
One has
\begin{align}
\chi(\pi'_k\tensor \pi'_{k+1})&=\chi(\pi_k\tensor s_j \pi_{k+1})=
\chi(\pi_k\tensor \pi_{k+1})-\<\alpha_j^\vee,\nu_{k+1}>\delta_{j,0}
\lbl{30a}\\
\chi(\pi'_r\tensor \pi'_{r+1})&=
\chi(s_j \pi_r\tensor s_j \pi_{r+1})\nonumber\\&=
\chi(\pi_r\tensor \pi_{r+1})+\<\alpha_j^\vee,\nu_{r}-\nu_{r+1}>
\delta_{j,0},
\qquad k<r<l\lbl{30b}\\
\chi(\pi'_l\tensor \pi'_{l+1})&=\chi(s_j \pi_l\tensor \pi_{l+1})+
\<\alpha_j^\vee,\nu_l>\delta_{j,0}\lbl{30c}
\end{align}
\end{clm}
Before we establish the claim, 
let us prove that~\loceqref{10a}--\loceqref{10c} follow from
\loceqref{30a}--\loceqref{30c} and~\loceqref{20}. 
Observe first that for~$j\not=0$, the above formulae imply that
$\chi(\pi'_r\tensor \pi'_{r+1})=
\chi(\pi_r\tensor \pi_{r+1})$ for all~$1\le r<Nm$, so
in that case there is nothing to prove. Furthermore, 
suppose that~$j=0$ and~$1\le l <Nm$. Then 
\begin{equation}\lbl{35}
\begin{split}
\Maj_\chi&(T_N(e_0 b))-\Maj_\chi(T_N(b))=
\sum_{r=k+1}^{l} r\<\alpha_0^\vee,\nu_{r}>
-\sum_{r=k}^{l-1}r\<\alpha_0^\vee,\nu_{r+1}>\\
&=\sum_{r=k+1}^l \<\alpha_0^\vee,\nu_r>=
N (h^0_b(k/Nm)-h^0_b(l/Nm))=-N
\end{split}
\end{equation}
by the choice of~$k$ and~$l$, whence by~\loceqref{30a}--\loceqref{30c}
\begin{equation}\lbl{40}
\begin{split}
K_s'-K_s=&\frac1N \sum_{r=k}^{\min\{s-1,l\}}
r(\chi(\pi_r\tensor \pi_{r+1})-\chi(\pi'_r\tensor \pi'_{r+1}))+\\
&\frac sN \sum_{r=\max\{s,k\}}^{l}
(\chi(\pi_r\tensor \pi_{r+1})-\chi(\pi'_r\tensor \pi'_{r+1})).
\end{split}
\end{equation}
For~$s=1,\dots,k$ the first sum is empty whilst the second sum reduces
to
$$
\sum_{r=k}^{l-1} \<\alpha_0^\vee,\nu_{r+1}>-
\sum_{r=k+1}^l \<\alpha_0^\vee,\nu_r>=0.
$$ 
Thus, $K_s'=K_s$, $s=1,\dots,k$. Furthermore, for~$s=k+1,\dots,l-1$ we get
\begin{align*}
K_s'-K_s=&\frac1N\Big(\sum_{r=k}^{s-1} r\<\alpha_0^\vee,
\nu_{r+1}>-\sum_{r=k+1}^s r\<\alpha_0^\vee,\nu_r>\Big)+
\frac sN\Big(\sum_{r=s}^{l-1}\<\alpha_0^\vee,\nu_{r+1}>-
\sum_{r=s}^l \<\alpha_0^\vee,\nu_r>\Big)\\
=&-\frac1N \sum_{r=k+1}^s \<\alpha_0^\vee,\nu_r>=
\<\alpha_0^\vee,\lambda_k-\lambda_s>.
\end{align*}
Finally, for~$s=l,\dots,Nm-1$ the second sum in~\loceqref{40} vanishes
and so
$$
K_s'-K_s=\frac1N \Big(\sum_{r=k}^{l-1} r \<\alpha_0^\vee,
\nu_{r+1}>-\sum_{r=k+1}^l r\<\alpha_0^\vee,\nu_r>\Big)=
-\frac1N \sum_{r=k+1}^l \<\alpha_0^\vee,\nu_r>=1
$$
by~\loceqref{35}.

Similarly, for~$l=Nm$ we get
$$
\Maj_\chi(T_N(e_0 b))-\Maj_\chi(T_N(b))=
\sum_{r=k+1}^l \alpha_0^\vee(\nu_r)-l\<\alpha_0^\vee,\nu_l>
=-N-Nm\<\alpha_0^\vee,\nu_l>,
$$
whence
\begin{equation}\lbl{50}
\begin{split}
K_s'-K_s=-\frac sN\,\<\alpha_0^\vee,\nu_{Nm}>+&\frac1N \sum_{r=k}^{s-1}
r(\chi(\pi_r\tensor \pi_{r+1})-\chi(\pi'_r\tensor \pi'_{r+1}))+\\
&\frac sN \sum_{r=\max\{s,k\}}^{Nm-1}
(\chi(\pi_r\tensor \pi_{r+1})-\chi(\pi'_r\tensor \pi'_{r+1})).
\end{split}
\end{equation}
For~$s=1,\dots,k$ the first sum in~\loceqref{50} is empty and so
$$
K_s'-K_s=-\frac sN\,\<\alpha_0^\vee,\nu_{Nm}>
+\frac sN \sum_{r=k}^{Nm-1} \<\alpha_0^\vee,\nu_{r+1}>-
\frac sN \sum_{r=k+1}^{Nm-1} \<\alpha_0^\vee,\nu_r>=0.
$$
Finally, for~$s=k+1,\dots,Nm$,
$$
K_s'-K_s=-\frac1N \sum_{r=k+1}^s \<\alpha_0^\vee,\nu_r>=
\<\alpha_0^\vee,\lambda_k-\lambda_s>,
$$
as required.

It remains to prove the claim. 
By Kashiwara's tensor product rule, $e_j b=b_1\tensor\cdots\tensor 
e_j b_p\tensor\cdots\tensor b_m$ for some~$1\le p\le m$.
Since~$b_p$ is an LS-path, the function~$h^j_{b_p}$
is strictly increasing on the interval~$[e^j_-(b_p),e^j_+(b_p)]=[k'/N,l'/N]$,
where~$k'=k-(p-1)N$, $l'=l-(p-1)N$ and~$0\le k'<l'<N$ 
by~\cite[Proposition~4.7(b)]{Li95} 
(recall that the definition of crystal operators on~$\hatbp$ we use 
differs by the sign of~$h^j$ from that
of~\cite{Li95}). Therefore, $\<\alpha_j^\vee,\nu_r><0$,
$r=k+1,\dots,l$. It follows from~\ref{PM40} that~$0=\varphi_j(\pi_r)<
\varepsilon_j(\pi_r)=-\alpha_j^\vee(\nu_r)$. Then,
by Kashiwara's tensor product rule
$$
e_j^{-\<\alpha_j^\vee,\nu_{r+1}>}(\pi_r\tensor
\pi_{{r+1}})=\pi_{r}\tensor e_j^{-\<\alpha_j^\vee,\nu_{r+1}>}
\pi_{{r+1}}=\pi_{r}\tensor s_j\pi_{{r+1}},\qquad
k<r<l,
$$
where we used~\ref{PM40}. Therefore, 
$\chi(\pi_{r}\tensor s_j\pi_{{r+1}})=\chi(\pi_{r}\tensor
\pi_{{r+1}})-\delta_{j,0}\<\alpha_j^\vee,\nu_{r+1}>$ by~\lemref{ZCB70}.
Furthermore, since~$\varepsilon_j(s_j\pi_{{r+1}})=\varepsilon_j(\pi_{r+1})
+\alpha_j^\vee(\nu_{r+1})=0$,
$$
e_j^{-\<\alpha_j^\vee,\nu_r>}(\pi_{r}\tensor s_j\pi_{{r+1}})=
e_j^{-\<\alpha_j^\vee,\nu_r>}\pi_{r}\tensor s_j\pi_{{r+1}}=
s_j\pi_{r}\tensor s_j\pi_{{r+1}}.
$$
Then~$\chi(s_j\pi_{r}\tensor s_j\pi_{{r+1}})=
\chi(\pi_{r}\tensor
s_j\pi_{{r+1}})+\delta_{j,0}\<\alpha_j^\vee,\nu_{r}>$
by~\lemref{ZCB70}, whence~\loceqref{30b}. 

Suppose that~$k>0$. Since all local maxima of~$h^j_b(\tau)$ are
integers, it follows by the choice of~$k$ and~$l$ that $h^j_{b}(\tau)
\le h^j_{b}(k/Nm)$ for~$\tau\le k/Nm$. 
Then~$\<\alpha_j^\vee,\nu_{k}>\le 0$, whence~$0=\varphi_j(\pi_{{k}})
<\varepsilon_j(\pi_{{k+1}})=-\<\alpha_j^\vee,\nu_{k+1}>$.
Using Kashiwara's tensor product rule, we obtain
$$
e_j^{-\<\alpha_j^\vee,\nu_{k+1}>}(\pi_{{k}}\tensor
\pi_{{k+1}})=\pi_{{k}} \tensor e_j^{-\<\alpha_j^\vee,\nu_{k+1}>}
\pi_{{k+1}}=\pi_{{k}} \tensor s_j\pi_{{k+1}},
$$
and~\loceqref{30a} follows by~\lemref{ZCB70}. Finally, suppose that~$l<Nm$.
Since~$h^j_{b}$ reaches its local maximum at~$l/Nm$,
we conclude that~$\<\alpha_j^\vee,\nu_{l+1}>\ge0$, whence
$\varphi_j(\pi_{{l}})=0=\varepsilon_j(\pi_{{l+1}})$. Therefore,
$$
e_j^{-\<\alpha_j^\vee,\nu_{l}>}(\pi_{{l}}\tensor \pi_{{l+1}})
=e_j^{-\<\alpha_j^\vee,\nu_{l}>}\pi_{{l}}\tensor \pi_{{l+1}}
=s_j \pi_{{l}}\tensor \pi_{{l+1}},
$$
which yields~\loceqref{30c} by~\lemref{ZCB70}.
\end{pf}
\begin{cor}
The associated crystal of a $\zeta$-crystal basis~$\wh B(\bsm\varpi_{i;m})$
of the integrable $\wh\bu_q$-module~$\wh V(\bsm\varpi_{i;m})$
is isomorphic to~$\psi(\wh{B(\varpi_i)^{\tensor m}})$.
\end{cor}

\subl{PM80}
Let us compute~$\psi(\pi_{\varpi_i}^{\tensor m}\tensor t^n)$, $n\in\bz$.
Recall that~$\chi(\pi_{\varpi_i}\tensor \pi_{\varpi_i})=0$. Since,
evidently, $T_N(\pi_{\varpi_i}^{\tensor m})=\pi_{\varpi_i}^{\tensor Nm}$
we conclude that~$\Maj_\chi(T_N(\pi_{\varpi_i}^{\tensor m}))=0$ and so
$\kappa_s(\pi_{\varpi_i}^{\tensor m}\tensor t^n)=s n/Nm$. Thus,
$\psi(\pi_{\varpi_i}^{\tensor m}\tensor t^n)=\pi_{m\varpi_i+n\delta}$.
In particular, $B(m\varpi_i+n\delta,\hatbp)$, $n\in\bz$ is an 
indecomposable subcrystal of~$\psi(\wh{B(\varpi_i)^{\tensor m}})$.
\begin{prop}
The image of~$\wh{B(\varpi_i)^{\tensor m}}$ in~$\hatbp$ 
under~$\psi$ is a disjoint
union of indecomposable crystals~$B(m\varpi_i+n\delta,\hatbp)$, 
$n=0,\dots,m-1$.
\end{prop}
\begin{pf}
By~\cite[Theorem~1.1]{NS}, a crystal basis of~$\wh V(\bsm\varpi_{i;1})$
is isomorphic to~$B(\varpi_i,\hatbp)$ which is indecomposable,
being generated by the linear path~$\pi_{\varpi_i}\in\hatbp$. 
On the other hand, 
$\wh B(\bsm\varpi_{i;1})$ is isomorphic to the
affinisation~$\wh{B(\bsm\varpi_{i;1})}$ of the crystal 
basis~$B(\bsm\varpi_{i;1})$ of~$V(\bsm\varpi_{i;1})$ which in turn
is isomorphic to the affinisation of~$B(\varpi_i)\subset \mbp$.
We conclude that~$\wh{B(\varpi_i)}$ is indecomposable as a crystal.
Thus there exists a monomial
$x\in \Cm$ such that~$x(\pi_{\varpi_i}\tensor 1)=\pi_{\varpi_i}\tensor t$.

Set~$\deg e_j=\delta_{j,0}$ and~$\deg f_j=-\delta_{j,0}$. That
defines a grading on~$\Cm$. Since~$x(\pi_{\varpi_i}\tensor 1)=
x\pi_{\varpi_i}\tensor t^{\deg x}$ by the definition of affinisation,
it follows that~$\deg x=1$ and~$x\pi_{\varpi_i}=\pi_{\varpi_i}$.
Write~$x=x_{j_k}\cdots x_{j_1}$ where~$x_{j_r}$ is either~$e_{j_r}$ 
or~$f_{j_r}$ and~$j_r\in\wh I$ and set~$x^{(m)}=x_{j_k}^m
\cdots x_{j_1}^m$. Evidently, $\deg x^{(m)}=m$.
We claim that~$x^{(m)} \pi_{\varpi_i}^{\tensor m}=\pi_{\varpi_i}^{\tensor m}$.
Indeed, $\pi_{\varpi_i}^{\tensor m}=S_m(\pi_{\varpi_i})$
and so~$x_{j_1}^m \pi_{\varpi_i}^{\tensor m}=S_m( x_{j_1} \pi_{\varpi_i})$.
It remains to use induction on~$k$.

Since~$\deg x^{(m)}=m$ it follows 
that~$x^{(m)}(\pi_{\varpi_i}^{\tensor m}\tensor t^k)=\pi_{\varpi_i}^{
\tensor m}\tensor t^{k+m}$ for all~$k\in\bz$.
On the other hand, let~$y=x'_{j_1}\cdots x'_{j_k}$ where~$x'_{j_r}=
e_{j_r}$ if~$x_{j_r}=f_{j_r}$ and vice versa. Evidently, $\deg y=-1$.
Since~$x\pi_{\varpi_i}=\pi_{\varpi_i}$ and~$e_i$, $f_i$ are
pseudo-inverses of each other, it follows that~$y\pi_{\varpi_i}=
\pi_{\varpi_i}$. Since~$\deg y^{(m)}=-m$ we conclude, as above,
that~$y^{(m)}(\pi_{\varpi_i}^{\tensor m}\tensor t^k)=\pi_{\varpi_i}^{
\tensor m}\tensor t^{k-m}$ for all~$k\in\bz$.
Therefore, $\pi_{m\varpi_i+(n+rm)\delta}
\in \hatbp$ lies in~$B(m\varpi_i+n\delta,\hatbp)$ for all~$r\in\bz$. 
It follows that~$B(m\varpi_i+r\delta,\hatbp)=B(m\varpi_i+s\delta,
\hatbp)$ if~$r=s\pmod m$.

Since~$B(\varpi_i)^{\tensor m}$ is indecomposable,
for all~$b\in B(\varpi_i)^{\tensor m}$ there exists a monomial~$x
\in \Cm$ such that~$b=x\pi_{\varpi_i}^{\tensor m}$.
Then for all~$k\in\bz$, $b\tensor t^k=x(\pi_{\varpi_i}^{\tensor m}
\tensor t^{k-\deg x})$. 
It follows that~$\psi(b\tensor t^k)\in B(m\varpi_i+n\delta,\hatbp)$
for some~$n\in\bz$, that is
$\psi(\wh{B(\varpi_i)^{\tensor m}})=\bigcup_{n=0}^{m-1}
B(m\varpi_i+n\delta,\hatbp)$. 

Since the 
crystals~$B(m\varpi_i+n\delta,\hatbp)$ are indecomposable, it remains
to prove that~$B(m\varpi_i+r\delta,\hatbp)\not= B(m\varpi_i+s\delta,
\hatbp)$, $r\not=s\pmod m$. For, observe that by the proof
of~\lemref{ZCB80} $\Maj_\chi(x b)=\Maj_\chi(b)-\deg x\pmod m$
provided that~$x b\not=0$, $x\in\Cm$.
Therefore, the set~$C_s:=
\{b\tensor t^k\,:\,b\in B(\varpi_i)^{\tensor m},\,
k\in\bz,\,\Maj_\chi(b)+k=s\pmod m\}$
is a subcrystal of~$\wh{B(\varpi_i)^{\tensor m}}$ and
$\wh{B(\varpi_i)^{\tensor m}}=\coprod_{s=0}^{m-1} C_s$.
Moreover, since~$\Maj_\chi(\pi_{\varpi_i}^{\tensor m})=0$,
$\pi_{\varpi_i}^{\tensor m}\tensor t^r\in C_s$ if and only if~$r=s
\pmod m$. Thus, $\psi(C_s)$ contains an element of~$B(m\varpi_i+r\delta,
\hatbp)$ if and only if~$r=s\pmod m$. It follows that~$\psi(C_s)=
B(m\varpi_i+s\delta)$.
\end{pf}
\begin{cor}[\thmref{thm2}]
The associated crystal of $\zeta$-crystal basis~$\wh B(\bsm\varpi_{i;m})^{(k)}$
of the simple $\wh\bu_q$-module $\wh V(\bsm\varpi_{i;m})^{(k)}$ is
isomorphic to~$B(m\varpi_i+k\delta,\hatbp)$ and hence indecomposable.
\end{cor}
\begin{pf}
By~\thmref{ZCB90}, $\wh B(\bsm\varpi_{i;m})^{(k)}=
\{b\tensor t^n\,:\, b\in B(\bsm\varpi_{i;m}),\,
n\in\bz,\,\Maj_\chi(b)+n=k\pmod m\}$. Since the associated crystal
of~$B(\bsm\varpi_{i;m})$ is isomorphic to~$B(\varpi_i)^{\tensor m}$
by~\cite{NS} and~\propref{ZCB30}, 
we conclude that~$\wh B(\bsm\varpi_{i;m})^{(k)}$
is isomorphic to~$C_k$ in the notation of the proof of the above
proposition. 
\end{pf}

\bibliographystyle{amsplain}

\end{document}